\documentclass[UKenglish,3p,sort&compress,final]{elsarticle}

\pdfoutput=1

\usepackage{babel}                    
 \usepackage{graphicx}
 \usepackage{amssymb,amsmath}
 \usepackage{hyperref}
 \usepackage{tikz}
\usepackage{pgfplots}
\pgfplotsset{compat=1.5}
 \usepackage{bbm}

 \usepackage{ifdraft}
 \usepackage{microtype}
 \usepackage{algpseudocode}
 \usepackage{mathtools}
 
\newcommand{\morespacearray}{\renewcommand{\arraystretch}{1.5}}
\newcommand{\g}{\ensuremath{\mathfrak{g}}}
\renewcommand{\d}{\ensuremath{\mathrm{d}}}
\DeclareMathOperator{\ad}{ad}
\DeclareMathOperator{\Ad}{Ad}
\DeclareMathOperator{\dexp}{dexp}
\DeclareMathOperator{\grad}{grad}

\DeclareMathOperator{\LieO}{O}
\DeclareMathOperator{\SO}{SO}
\DeclareMathOperator{\SL}{SL}
\DeclareMathOperator{\SE}{SE}
\DeclareMathOperator{\GL}{GL}
\DeclareMathOperator{\Sp}{Sp}

\DeclareMathOperator{\St}{St}

\newcommand{\Id}{\mathrm{Id}}

\newcommand{\R}{\mathbf{R}}
\newcommand{\mdef}{\coloneqq}
\newcommand{\fedm}{\eqqcolon}
\newcommand{\dd}{\ensuremath{\overline{\d}}}
\newcommand{\q}{\mathbbm{q}}

\journal{Journal of Computational Physics}
  \begin{document}
  
 \begin{frontmatter} 
 \title{An introduction to Lie group integrators -- basics, new developments and applications}
 \author[ntnu]{Elena Celledoni}
 \ead{elenac@math.ntnu.no}
 \author[ntnu]{H{\aa}kon Marthinsen}
\ead{hakonm@math.ntnu.no}
 \author[ntnu]{Brynjulf Owren}
 \ead{bryn@math.ntnu.no}  
 
 \address[ntnu]{Department of Mathematical Sciences, NTNU, N--7491 Trondheim, Norway}

\begin{abstract}
We give a short and elementary introduction to Lie group methods. 
A selection of applications of Lie group integrators are discussed.
Finally, a family of symplectic integrators on cotangent bundles of Lie groups is presented and the notion of discrete gradient methods is generalised to Lie groups.
\end{abstract}

\end{frontmatter}

\section{Introduction}
The significance of the geometry of differential equations was well understood already in the nineteenth century, and in the last few decades such aspects have played an increasing role in numerical methods for differential equations. Nowadays, there is a rich selection of integrators which preserve properties like symplecticity, reversibility, phase volume and first integrals, either exactly or approximately over long times~\cite{hairer10gni}. Differential equations are inherently connected to Lie groups, and in fact one often sees applications in which the phase space is a Lie group or a manifold with a Lie group action. In the early nineties, two important papers appeared which used the Lie group structure directly as a building block in the numerical methods. Crouch and Grossman \cite{crouch93nio} suggested to advance the numerical solution by computing flows of vector fields in some Lie algebra.  Lewis and Simo \cite{lewis94caf} wrote an influential paper on Lie group based integrators for Hamiltonian problems, considering the preservation of symplecticity, momentum and energy. These ideas were developed in a systematic way throughout the nineties by several authors. In a series of three papers, Munthe-Kaas \cite{munthe-kaas95lbt, munthe-kaas98rkm,munthe-kaas99hor} presented what are now known as the Runge--Kutta--Munthe-Kaas methods. By the turn of the millennium, a survey paper \cite{iserles00lgm} summarised most of what was known by then about Lie group integrators. More recently a survey paper on structure preservation appeared with part of it dedicated to Lie group integrators \cite{christiansen11tis}.

 The purpose of the present paper is three-fold. First, in section~\ref{sec:basics} we give an elementary, geometric introduction to the ideas behind Lie group integrators. Secondly, we present some examples of applications of Lie group integrators in sections
\ref{sec:appcm} and \ref{sec:dataanalysis}. There are many such examples to choose from, and we give here only a few teasers.
These first four sections should be read as a survey.
But in the last two section, new material is presented.
Symplectic Lie group integrators have been known for some time, derived by Marsden and coauthors \cite{marsden01dma}
by means of variational principles. In section~\ref{sec:symplie} we consider a group structure on the cotangent bundle of a Lie group and derive symplectic Lie group integrators using the model for vector fields on manifolds defined by Munthe-Kaas in~\cite{munthe-kaas99hor}.
In section~\ref{sec:discdiff} we extend the notion of discrete gradient methods as proposed by
Gonzalez \cite{gonzalez96tia} to Lie groups, and thereby we obtain a general method for preserving first integrals in differential equations on Lie groups.

We would also like to briefly mention some of the issues we are \emph{not} pursuing in this article. One is the important family of Lie group integrators for problems of linear type, including methods based on the Magnus and Fer expansions. An excellent review of the history, theory and applications of such integrators can be found in \cite{blanes09tme}.
We will also skip all discussions of order analysis of Lie group integrators.
This is a large area by itself which involves technical tools and mathematical theory which we do not wish to include in this relatively elementary exposition.
There have been several new developments in this area recently, in particular by Lundervold and Munthe-Kaas, see e.g.\ \cite{lundervold11hao}.

\section{Lie group integrators}

\label{sec:basics}
The simplest consistent method for solving ordinary differential equations is the Euler method.
For an initial value problem of the form
\begin{equation*}
   \dot{y} = F(y),\quad y(0)=y_0,
\end{equation*}
 one takes a small time increment $h$, and approximates $y(h)$ by the simple formula
 \begin{equation*}
     y_{1} = y_0 + hF(y_0),
 \end{equation*}
advancing along the straight line coinciding with the tangent at $y_0$.
Another way of thinking about the Euler method is to consider the constant vector field
$F_{y_0}(y) \mdef F(y_0)$ obtained by parallel translating the vector~$F(y_0)$ to all points  of phase space.
A step of the Euler method is nothing else than computing the exact $h$-flow of this simple vector field starting at $y_0$.  In Lie group integrators, the  same principle is used, but allowing for more advanced vector fields than the constant ones. A Lie group generalisation of the Euler method is called the Lie--Euler method, and we shall illustrate its use through an example \cite{crouch93nio}.
\paragraph{Example, the Duffing equation}
Consider the system in $\R^2$
\begin{equation} \label{eq:duffing}
\begin{aligned}
    \dot{x} &= y, \\
    \dot{y} &= -a x - b x^3, 
\end{aligned} \qquad a\geq 0, b\geq 0,
\end{equation}
a model used to describe the buckling of an elastic beam. Locally, near a point
$(x_0,y_0)$ we could use the approximate system
\begin{equation}
\begin{alignedat}{2} \label{eq:duff:sl2:frozen}
   \dot{x} &= y,              &\qquad x(0) &= x_0,\\
   \dot{y} &= -(a+b x_0^2) x, &       y(0) &= y_0,
\end{alignedat} 
\end{equation}
which has the exact solution
\begin{equation}\label{eq:duffing:sl2:flow}
      \bar{x}(t) = x_0\cos\omega t+ \frac{y_0}{\omega}\sin\omega t,\quad
       \bar{y}(t) = y_0\cos\omega t- \omega x_0\sin\omega t,\quad \omega=\sqrt{a+bx_0^2}.
\end{equation}
Alternatively, we may consider
   the local problem
\[
\begin{aligned}
   \dot{x} &= y, \\
   \dot{y} &= -ax  - bx_0^3,
\end{aligned}
\]
having exact solution
\[
\begin{aligned}
\bar{x}(t) &= x_0\cos\alpha t + \frac{y_0}{\alpha}\sin\alpha t + b\,x_0^3\,\frac{\cos\alpha t-1}{\alpha^2},\\
\bar{y}(t) &= y_0\cos\alpha t - \alpha x_0\sin\alpha t - b\,x_0^3\, \frac{\sin\alpha t}{\alpha},
\end{aligned} \qquad \alpha=\sqrt{a}.
\]
In each of the two cases, one may take $x_1=\bar{x}(h)$, $y_1=\bar{y}(h)$ as the numerical approximation at time $t=h$. 
The same procedure is repeated in subsequent steps.
\begin{figure}
\centering
\includegraphics{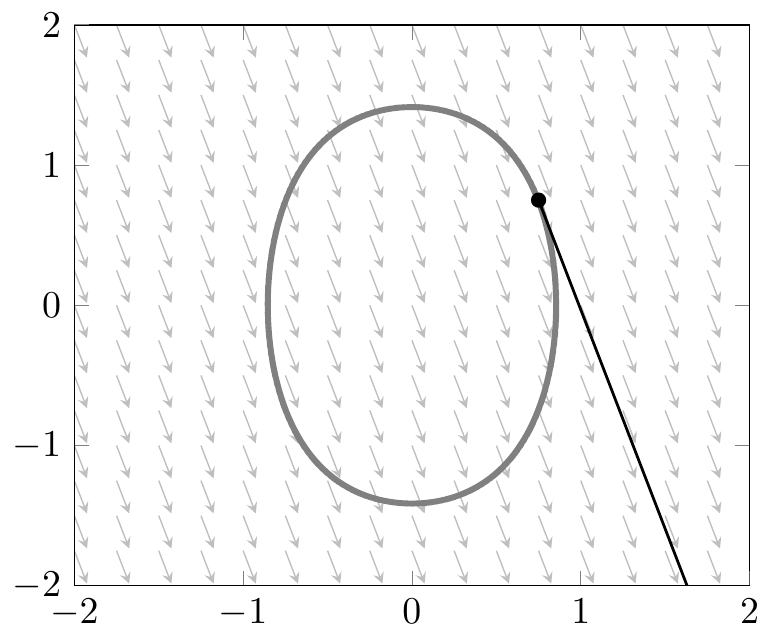}
\includegraphics{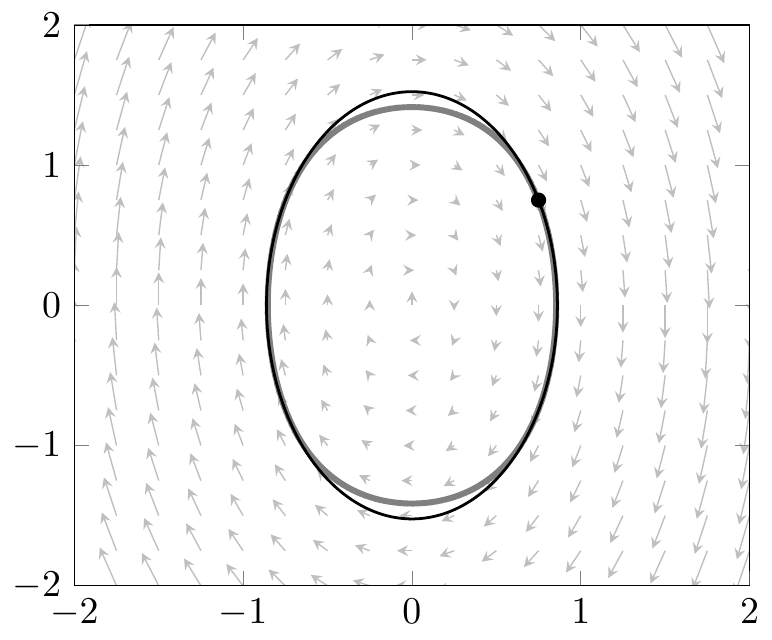}
\caption{$(\R^d,+)$-frozen vector field (left) and $\mathfrak{sl}(2)$-frozen vector field (right) for the Duffing equation. Both are frozen at
$(x_0,y_0)=(0.75,0.75)$.
The thin black curve in each plot shows the flows of the frozen vector fields $0\leq t\leq 20$. The thicker curve in each plot is the exact flow of the Duffing equation.
 \label{fig:duffing}}
 \end{figure}
 A common framework for discussing these two cases is provided by the use of frames, i.e.\ a set of 
 of vector fields which at each point is spanning the tangent space.
 In the first case, the numerical method applies the frame 
\begin{equation}  \label{eq:XYsl2}
      X = \begin{bmatrix} y\\0 \end{bmatrix} \fedm y\, \partial x,\quad
      Y = \begin{bmatrix} 0\\x \end{bmatrix} \fedm x\, \partial y.
\end{equation}
Taking the negative Jacobi--Lie bracket (also called the commutator) between $X$ and $Y$ yields the third element of the standard basis for the Lie algebra $\mathfrak{sl}(2)$, i.e.\
\begin{equation} \label{eq:Hsl2}
    H = -[X,Y] =  x\,\partial x - y\,\partial y,
\end{equation}
so that the frame may be augmented to consist of $\{X, Y, H\}$.
 In the second case, the vector fields $E_1=y\,\partial x - ax\,\partial y$ and $E_2=\partial y$ can  be used as
a frame, but again we choose to  augment these two fields with the commutator $E_3=-[E_1,E_2]=\partial x$ to obtain the Lie algebra of the special Euclidean group $\SE(2)$ consisting of translations and rotations in the plane.  The situation is illustrated in Figure~\ref{fig:duffing}.
In the left part, we have considered the constant vector field corresponding to the Duffing vector field evaluated at
$(x_0,y_0)=(0.75,0.75)$, and the exact flow of this constant field is just the usual Euler method, a straight line.
In the right part, we have  plotted the vector field defined in \eqref{eq:duff:sl2:frozen} with the same $(x_0,y_0)$ along with its
flow~\eqref{eq:duffing:sl2:flow}. The exact flow of \eqref{eq:duffing} is shown in both plots (thick curve). 
 
In general, a way to think about Lie group integrators is that we have a manifold $M$ where there is such a frame available; $\{E_1,\ldots,E_d\}$ such that at any point $p\in M$ one has
\[
   \operatorname{span}\{E_1(p),\ldots,E_d(p)\} = T_pM.
\]
Frames with this property are said to be locally transitive.
The frame may be a linear space or in many cases even a Lie algebra $\g$ of vector fields.
In the example with Duffing's equation, the set $\{X, Y, H\}$ is locally transitive on $\R^2\setminus\{0\}$ and $\{E_1,E_2,E_3\}$ is locally transitive on $\R^2$.

Given an arbitrary vector field $F$ on $M$, then at any point
$p\in M$ there exists a vector field $F_p$ in the span of the frame vector fields
such that $F_p(p) = F(p)$.  
An explicit way of writing this is by using a set of basis vector fields $E_1,\ldots,E_d$ for $\g$, such that any smooth vector field
$F$ has a representation
\begin{equation} \label{VFframerep}
     F(y) = \sum_{k=1}^d f_k(y) E_k(y),
\end{equation}
for some functions $f_k \colon M\rightarrow\R$. The vector fields $F_p\in\g$, called 
\emph{vector fields with frozen coefficients} by Crouch and Grossman \cite{crouch93nio}, are then obtained
as
\[
     F_p(y) = \sum_{k=1}^d f_k(p) E_k(y).
\]
In the example with the Duffing equation we took $E_1=X, E_2=Y$,
 $f_1(x,y)=1$ and $f_2(x,y)=-(a+bx^2)$. The Lie--Euler method reads in general
\begin{equation} \label{eq:Lie--Euler}
      y_{n+1} = \exp(hF_{y_n})y_n,
\end{equation}
where $\exp$ denotes the flow of a vector field.

A more interesting example, also found in \cite{crouch93nio} is obtained by choosing $M=S^2$, the 2-sphere.
A suitable way to induce movements of the sphere is that of rotations, that is, by introducing the Lie group $\SO(3)$ consisting of orthogonal matrices with unit determinant. The corresponding Lie algebra $\mathfrak{so}(3)$ of vector fields are spanned by
\[
      E_1(x,y,z)=-z\,\partial y + y\,\partial z,\quad
      E_2(x,y,z)=z\,\partial x - x\,\partial z,\quad
      E_3(x,y,z)=-y\,\partial x + x\,\partial y.
\]
We note that $xE_1(x,y,z)+yE_2(x,y,z)+zE_3(x,y,z)=0$, showing that the functions $f_k$
in \eqref{VFframerep} are not unique.
A famous example of a system whose solution evolves on $S^2$ is the free rigid body Euler equations
\begin{equation} \label{eq:euler_frb}
    \dot{x} = \left(\frac1{I_3}-\frac1{I_2}\right)\,yz,  \quad
    \dot{y} = \left(\frac1{I_1}-\frac1{I_3}\right)\,xz, \quad
  \dot{z} = \left(\frac1{I_2}-\frac1{I_1}\right)\,xy,
\end{equation}
where $x, y, z$ are the coordinates of the angular momentum relative to the body, and $I_1, I_2, I_3$ are the principal moments of inertia. A choice of representation \eqref{VFframerep} is obtained with
\[
    f_1(x,y,z) = -\frac{x}{I_1},\quad f_2(x,y,z)=-\frac{y}{I_2},\quad f_3(x,y,z) = -\frac{z}{I_3},
\]
so that the ODE vector field can be expressed in the form
\[
F(x,y,z)=
-\frac{x}{I_1}
\begin{bmatrix*}[r]
0\\ -z\\ y
\end{bmatrix*}
-\frac{y}{I_2}
\begin{bmatrix*}[r]
z\\ 0\\ -x
\end{bmatrix*}
-\frac{z}{I_3}
\begin{bmatrix*}[r]
-y\\ x\\ 0
\end{bmatrix*}.
\]
We compute the vector field with coefficients frozen at $p_0=(x_0,y_0,z_0)$,
\[
F_{p_0}(x,y,z)= \mathbf{F}_{p_0} \begin{bmatrix} x\\ y \\ z \end{bmatrix}
\mdef \begin{bmatrix*}[r]
0 & \frac{z_0}{I_3} & -\frac{y_0}{I_2}  \\
-\frac{z_0}{I_3} & 0 & \frac{x_0}{I_1} \\
\frac{y_0}{I_2} & -\frac{x_0}{I_1} & 0
\end{bmatrix*}
\begin{bmatrix} x\\ y \\ z \end{bmatrix}.
\]
The $h$-flow of this vector field is the solution of a linear system of ODEs and can be expressed in terms of the matrix exponential $\mathtt{expm}(h\mathbf{F}_{p_0})$. The Lie--Euler method can be expressed
as follows: 
\medskip
\begin{algorithmic}
	\State $p_0 \gets (x_0,y_0,z_0)$
	\For{$n \gets 0, 1, \dotsc$}
		\State $p_{n+1} \gets \mathtt{expm}(h\mathbf{F}_{p_n}) p_n$
	\EndFor
\end{algorithmic}
\medskip
Notice that the matrix to be exponentiated belongs to the matrix group $\mathfrak{so}(3)$ of real skew-symmetric matrices. The celebrated Rodrigues' formula
\[
  \mathtt{expm}(A) = I + \frac{\sin\alpha}{\alpha} A + \frac{1-\cos\alpha}{\alpha^2}A^2,\qquad
  \alpha^2=\lVert A\rVert_2^2=\tfrac12\lVert A\rVert_F^2,\quad A\in\mathfrak{so}(3),
\]
provides an inexpensive way to compute this.

Whereas the notion of frames was used by Crouch and Grossman in their pioneering work
\cite{crouch93nio}, a different type of notation was used in a series of papers by Munthe-Kaas
\cite{munthe-kaas95lbt,munthe-kaas98rkm,munthe-kaas99hor}, see also \cite{lundervold11hao} for a more modern treatment. Let $G$ be a finite dimensional Lie group acting transitively on a manifold $M$. A Lie group action is generally a map from $G\times M$ into $M$, having the properties that
\begin{equation*} \label{eq:LGaction}
      e\cdot m=m,\;\forall m\in M,\qquad g\cdot(h\cdot m)=(g\cdot h)\cdot m,\ \forall g,h\in G,\; m\in M,
\end{equation*}
where $e$ is the group identity element, and  the first $\cdot$ in the right hand side of the second identity is the group product.
Transitivity means that for any two points $m_1, m_2\in M$ there exists a group element $g\in G$ such that $m_2=g\cdot m_1$. We denote the Lie algebra of $G$ by $\g$.
For any element $\xi\in\g$ there exists a vector field on~$M$
\begin{equation} \label{eq:groupaction}
   X_{\xi}(m) = \left.\frac{\d}{\d t}\right\rvert_{t=0} \exp(t\xi)\cdot m \fedm \lambda_*(\xi)(m).
\end{equation}
Munthe-Kaas introduced a generic representation of a vector field  $F\in\mathcal{X}(M)$   
by  a map $f \colon M\rightarrow\g$ such that
\begin{equation} \label{eq:genpres}
       F(m) =\lambda_*(f(m))(m).
\end{equation}
The corresponding frame is obtained as $E_i=\lambda_*(e_i)$  where $\{e_1,\ldots,e_d\}$ is some basis for $\g$ and one chooses the functions $f_i \colon M\rightarrow\R$ such that 
$f(m) = \sum_{i=1}^d f_i(m) e_i$. The map $\lambda_*$ is an anti-homomorphism of the Lie algebra $\g$ into the Lie algebra of vector fields $\mathcal{X}(M)$ under the Jacobi--Lie bracket, meaning that
\[
     \lambda_*([X_m, Y_m]_{\g}) =- [\lambda_*(X_m), \lambda_*(Y_m)]_{\mathrm{JL}}.
\]
This separation of the Lie algebra $\g$ from the manifold $M$ allows for more flexibility in the way we represent the frame vector fields. For instance, in the example with Duffing's equation and
the use of $\mathfrak{sl}(2)$, we could have used the matrix Lie algebra with basis elements
\[
     X_m =
     \begin{bmatrix}
     0 & 1 \\ 0 & 0
     \end{bmatrix},\quad
          Y_m =
     \begin{bmatrix}
     0 & 0 \\ 1 & 0
     \end{bmatrix},\quad
     H_m =
     \begin{bmatrix*}[r]
     1 & 0 \\ 0 & -1
     \end{bmatrix*},
\]
rather than the basis of vector fields \eqref{eq:XYsl2}, \eqref{eq:Hsl2}. The group action 
by $g\in \SL(2)$ on a point $m\in\R^2$ would then be simply $g\cdot m$, matrix-vector multiplication, and the $\exp$ in \eqref{eq:groupaction} would be the matrix exponential.
The map~$f(x,y)$ would in this case be
\[
      f \colon (x,y) \mapsto
\begin{bmatrix}
   0 & y \\ -(a+bx^2) & 0
  \end{bmatrix},
\]
but note that since the dimension of the manifold is just two whereas the dimension of
$\mathfrak{sl}(2)$ is three, there is freedom in the choice of $f$. In the example we chose not to use the third basis element $H$.

\subsection{Generalising Runge--Kutta methods}
\label{GeneralizingRK}
In order to construct  general schemes, as for instance a Lie group version of the Runge--Kutta methods, one needs to introduce intermediate stage values. This can be achieved in a number of different ways. They all have in common that when the methods are applied in Euclidean space where the Lie group is $(\R^m,+)$, they reduce to conventional Runge--Kutta schemes.
Let us begin by studying the simple second order Heun method, sometimes called the improved Euler method. 
\[
    k_1 = F(y_n),\quad k_2=F(y_n+hk_1),\qquad y_{n+1} = y_n + \tfrac12 h(k_1+k_2).
\]
 Geometrically,  we may think of $k_1$ and $k_2$ as constant vector fields, coinciding with the exact ODE $F(y)$ at the points $y_n$ and $y_n+hk_1$ respectively.
The update $y_{n+1}$ can be interpreted in at least three different ways,
\begin{equation} \label{eq:flowalts}
      \exp\left(\frac{h}{2}(k_1+k_2)\right) \cdot y_n, \quad
     \exp\left(\frac{h}{2}k_1\right) \cdot    \exp\left(\frac{h}{2}k_2\right) \cdot y_n, \quad 
         \exp\left(\frac{h}{2}k_2\right) \cdot    \exp\left(\frac{h}{2}k_1\right) \cdot y_n.
\end{equation}
 The first is an example of a Runge--Kutta--Munthe-Kaas method and the second is an example of a Crouch--Grossman method. All three fit into the framework of \emph{commutator-free} Lie group methods. All three suggestions above are generalisations  that will reduce to Heun's method in $(\R^m,+)$.
In principle we could extend the idea to Runge--Kutta methods with several stages
 \[
 y_{n+1} = y_n + h\sum_{i=1}^sb_iF(Y_i),\quad Y_i = y_n + h\sum_{j=1}^s a_{ij} F(Y_j),\
 i=1,\ldots,s,
 \]
  by for instance interpreting the summed expressions as vector fields with frozen coefficients whose flows we apply to the point $y_n\in M$.
 But it is unfortunately not true that one in this way will retain the order of the Runge--Kutta method when applied to cases where the acting group is non-abelian.

 Let us first describe methods as proposed by Munthe-Kaas \cite{munthe-kaas99hor}, where one may think of the method simply as a change of variable. As before, we assume that the action of $G$ on $M$ is locally transitive. Since the exponential mapping is a local diffeomorphism in some open set containing $0\in\g$, it is possible to represent \emph{any} smooth curve $y(t)$ on $M$ in some neighbourhood of a point $p\in M$ by means of a curve 
 $\sigma(t)$ through the origin of $\g$ as follows
 \begin{equation} \label{eq:solrep}
      y(t) = \exp(\sigma(t))\cdot p,\qquad \sigma(0)=0,
 \end{equation}
  though $\sigma(t)$ is not necessarily unique. We may differentiate this curve with respect to $t$ to obtain  
\begin{equation} \label{eq:doty}
       \dot{y}(t) =   \lambda_*\bigl(\dexp_{\sigma(t)}\dot{\sigma}(t)\bigr)(y(t)) = F(y(t))=
       \lambda_*\bigl(f(\exp(\sigma(t))\cdot p)\bigr)(y(t)).
\end{equation}
The details are given in \cite{munthe-kaas99hor} and the map 
$\dexp_\sigma \colon \g\rightarrow\g$ was derived by Hausdorff in \cite{hausdorff06dse} as an infinite series of commutators
\begin{equation} \label{eq:dexp}
    \dexp_{\sigma}(v) = v + \frac{1}{2}[\sigma,v] + \frac{1}{6}[\sigma,[\sigma,v]]+\dotsb
    = \sum_{k=0}^\infty \frac{1}{(k+1)!} \ad_\sigma^k v = \left.\frac{\exp(z)-1}{z}\right\rvert_{z=\ad_\sigma} v,
\end{equation}
with the usual definition of $\ad_u(v)$ as the commutator $[u,v]$.
The map $\lambda_*$ does not have to be injective, but a sufficient condition for \eqref{eq:doty} to hold is
that
\[
     \dot{\sigma} = \dexp_{\sigma}^{-1} (f(\exp(\sigma)\cdot p)).
\]  
This is a differential equation for $\sigma(t)$ on a linear space, and one may choose any conventional integrator for solving it. The map $\dexp_\sigma^{-1} \colon \g\rightarrow\g$ is the inverse of $\dexp_\sigma$ and can also be obtained by differentiating the logarithm, i.e.\ the inverse of the exponential map.
  From \eqref{eq:dexp} we find that one can write $\dexp_{\sigma}^{-1}(v)$ as
\begin{equation} \label{eq:dexpinv}
\dexp_{\sigma}^{-1}(v) = \left.\frac{z}{\exp(z)-1}\right\rvert_{z=\ad_\sigma}v
= v - \frac{1}{2}[\sigma,v] + \frac{1}{12}[\sigma, [\sigma, v]]+\cdots.
\end{equation}
The coefficients appearing in this expansion are scaled Bernoulli numbers $\frac{B_k}{k!}$, and  $B_{2k+1}=0$ for all $k\geq 1$. One step of the resulting Runge--Kutta--Munthe-Kaas method 
is then expressed in terms of evaluations of the map $f$ as follows
\begin{align*}
y_{1} &= \exp\Bigl(h\sum_{i=1}^s b_i k_i\Bigr) \cdot y_0, \\
k_i &= \dexp_{h\sum_j a_{ij}k_j}^{-1} f\biggl(\exp\Bigl(h\sum_j a_{ij}k_j\Bigr) \cdot y_0\biggr),\quad i=1,\ldots,s.
\end{align*}
This is not so surprising seen from the perspective of the first alternative in \eqref{eq:flowalts}, the main difference is that the stages $k_i$ corresponding to the frozen vector fields
$\lambda_*(k_i)$ need to be ``corrected'' by the $\dexp^{-1}$ map.
Including this map in computational algorithms may seem awkward, however, fortunately truncated versions of \eqref{eq:dexpinv} may frequently be used. In fact, by applying some clever tricks involving graded free Lie algebras, one can in many cases replace $\dexp^{-1}$ with a low order Lie polynomial while retaining the convergence order of the original Runge--Kutta method. Details of this can be found in \cite{munthe-kaas99cia,
casas03cel}.
 
 There are also some important cases of Lie algebras for which $\dexp_\sigma^{-1}$ can be computed exactly in terms of elementary functions, among those is $\mathfrak{so}(3)$ reported in \cite{celledoni03lgm}. Notice that the representation \eqref{eq:solrep} does not depend on the use of the exponential map from $\g$ to $G$. In principle, one can replace this map with any local diffeomorphism $\varphi$, where one usually scales $\varphi$ such that $\varphi(0)=e$ and
$T_0\varphi = \Id_\g$.  An example of such map is the Cayley transformation \cite{diele98tct} which can be used for matrix Lie groups of the type
$G_P = \{X\in\R^{n\times n} \mid X^TPX = P\}$ for a nonsingular $n\times n$-matrix $P$.
 These include the orthogonal group~$\LieO(n)=G_I$ and
the linear symplectic group~$\Sp(n)=G_J$ where $J$ the skew-symmetric matrix of the standard symplectic form. Another possibility is to replace the exponential map by canonical coordinates of the second kind \cite{owren01imb}.

We present here the well-known Runge--Kutta--Munthe-Kaas method based on the popular fourth order method of Kutta \cite{kutta01bzn}, having Butcher tableau
\begin{equation} \label{eq:kutta}
\morespacearray
\begin{array}{r|rrrr}
0 & \\
\tfrac{1}{2} & \tfrac{1}{2}  \\
\tfrac{1}{2} & 0 & \tfrac{1}{2} \\
1 & 0 & 0 & 1 \\ \hline
 & \tfrac{1}{6} &  \tfrac{1}{3} &  \tfrac{1}{3} &  \tfrac{1}{6}
\end{array}
\end{equation}
In the Lie group method, the 
$\dexp^{-1}$ map has been replaced by the optimal Lie polynomials.
\begin{align*}
   k_1 &= h f(y_0),\\
   k_2 &= h f(\exp(\tfrac{1}{2}k_1) \cdot y_0), \\
   k_3 &= h f(\exp(\tfrac{1}{2}k_2-\tfrac{1}{8}[k_1,k_2])\cdot y_0), \\
   k_4 &= h f(\exp(k_3)\cdot y_0), &  \\
   y_1 &= \exp(\tfrac{1}{6}(k_1+2k_2+2k_3+k_4-\tfrac12[k_1,k_4]))\cdot y_0.
\end{align*}
An important advantage of the Runge--Kutta--Munthe-Kaas schemes is that it is easy to preserve the convergence order when extending them to Lie group integrators.
This is not the case with for instance the schemes of Crouch and Grossman \cite{crouch93nio, owren99rkm}, where it is necessary to develop order conditions for the non-abelian case.
This is also true for the commutator-free methods developed by Celledoni et al.\ \cite{celledoni03cfl}.
In fact, these methods include those of Crouch and Grossman.
The idea here is to allow compositions of exponentials or flows instead of commutator corrections. With stages $k_1,\ldots, k_s$ in the Lie algebra, one includes expressions of the form
\[
      \exp\Bigl(\sum_{i}\beta_J^i k_i\Bigr) \dotsm \exp\Bigl(\sum_{i}\beta_2^i k_i\Bigr) \cdot \exp\Bigl(\sum_{i}\beta_1^i k_i\Bigr) \cdot y_0,
\]
both in the definition of the stages and the update itself.
In some cases it is also possible to reuse flow calculations from one stage to another, and thereby lower the computational cost of the scheme. An extension of \eqref{eq:kutta} can be obtained as follows, setting $k_i=hf(Y_i)$ for all $i$,
\begin{align*}
 Y_1&= y_0, \\ 
 Y_2&=\exp(\tfrac{1}{2}k_1)\cdot y_0, \\
 Y_3&=\exp(\tfrac{1}{2}k_2)\cdot y_0 \\
 Y_4&= \exp(k_3-\tfrac{1}{2}k_1)\cdot Y_2, \\
 y_{\frac{1}{2}} &= 
 \exp(\tfrac{1}{12}(3k_1+2k_2+2k_3-k_4))\cdot y_0, \\
 y_{1} &= 
 \exp(\tfrac{1}{12}(-k_1+2k_2+2k_3+3k_4))\cdot y_{\frac{1}{2}}. 
  \end{align*}
Note in particular in this example how the expression for $Y_4$ involves
$Y_2$ and thereby one exponential calculation has been saved.

\subsection{A plenitude of group actions}
We saw in the first examples with Duffing's equation that the manifold $M$, the group $G$ and even the way $G$ acts on $M$ can be chosen in different ways. It is not obvious which action is the best or suits the purpose in the problem at hand. Most examples we know from the literature are using matrix Lie groups~$G\subseteq \GL(n)$, but the choice of group action depends on the problem and the objectives of the simulation. We give here several examples of situations where Lie group integrators can be used.

\paragraph{$G$ acting on $G$} In the case $M=G$, it is natural to use either left or right multiplication as the action
\[
      L_g(m) = g\cdot m\quad\text{or}\quad R_g(m) = m\cdot g,\qquad g,m\in G.
\]
The correspondence between the vector field $F\in\mathcal{X}(M)$ and the
map \eqref{eq:genpres} is then just the tangent map of left or right multiplication
\[
     F(g) = T_eL_g(f(g))\quad\text{or}\quad F(g)=T_eR_{g}(\tilde{f}(g)),\quad g\in G.
\]
When working with matrices, this simply amounts to
setting $F(g) = g\cdot f(g)$ or $F(g) = \tilde{f}(g) \cdot g$. Note that $\tilde{f}(g)$ is related to 
$f(g)$ through the adjoint representation of $G$, $\Ad \colon G\rightarrow \operatorname{Aut}(\g)$,
 \[
      \tilde{f}(g) = \Ad_g f(g),\qquad \Ad_g = T_eL_g\circ T_eR_g^{-1}.
\]

\paragraph{The affine group and its use in semilinear PDE methods}
Lie group integrators can also be used for approximating the solution to partial differential equations, the most obvious choice of PDE model being the semilinear problem
\begin{equation} \label{eq:semilinear}
     u_t = Lu + N(u),
\end{equation}
where $L$ is a linear differential operator and $N(u)$ is some nonlinear map, typically containing derivatives of lower order than $L$. After discretising in space, \eqref{eq:semilinear} is turned into a system of  $n_d$ ODEs, for some large $n_d$, $L$ becomes an $n_d\times n_d$-matrix, and $N \colon \R^{n_d}\rightarrow \R^{n_d}$ a nonlinear function. 
We may now as in \cite{munthe-kaas99hor} introduce the action on $\R^{n_d}$ by some subgroup of the affine group represented as the semidirect product $G=\GL(n_d)\ltimes\R^{n_d}$.
The group product, identity, and inverse are given as
\[
   (A_1,b_1)\cdot (A_2,b_2) = (A_1 A_2, A_1b_2+b_1),\quad
   e=(I,0),\quad (A,b)^{-1}=(A^{-1}, -A^{-1}b).
\]
The action on $\R^{n_d}$ is
\[
     (A,b)\cdot x = Ax + b,\qquad (A,b)\in G,\ x\in\R^{n_d},
\]
and the Lie algebra and commutator are given as
\[
   \g = (\xi, c),\ \xi\in\mathfrak{gl}(n_d),\ c\in\R^{n_d},\quad
   [(\xi_1,c_1),(\xi_2,c_2)] = ([\xi_1,\xi_2], \xi_1c_2-\xi_2c_1+c_1).
\]
In many interesting PDEs, the operator $L$ is constant, so it makes sense
to consider the $n_d+1$-dimensional subalgebra $\g_L$ of $\g$ consisting of elements 
$(\alpha L, c)$ where $\alpha\in\R$, $c\in\R^d$, so that the
map $f \colon \R^{n_d}\rightarrow \g_L$ is given as
\[
      f(u) = (L,N(u)).
\]
One parameter subgroups are obtained through the exponential map as follows
\[
   \exp(t(L,b)) = (\exp(tL), \phi(tL) tb).
\]
Here the entire function $\phi(z)=(\exp(z)-1)/z$ familiar from the theory of exponential integrators appears. As an example, one could now consider the Lie--Euler method \eqref{eq:Lie--Euler} in this setting, which coincides with the exponential Euler method
\[
   u_{1} = \exp(h(L,N(u_0))\cdot u_0 = \exp(hL)u_0 + h\phi(hL)N(u_0).
\]
There is a large body of literature on exponential integrators, going almost half a century back in time, see~\cite{hochbruck10ei} and the references therein for an extensive account.

\paragraph{The coadjoint action and Lie--Poisson systems}
Lie group integrators for this interesting case were studied by Eng{\o} and Faltinsen \cite{engo01nio}.
Suppose $G$ is a Lie group and the manifold under consideration is the dual space~$\g^*$ of its Lie algebra~$\g$. The coadjoint action by $G$ on $\g^*$ is denoted 
$\Ad_g^*$ defined for any $g\in G$ as
\begin{equation} \label{eq:coadjoint-action}
\langle \Ad_g^*\mu, \xi\rangle = \langle \mu, \Ad_g\xi\rangle,\quad\forall\xi\in\g,
\end{equation}
for a duality pairing $\langle {\cdot}, {\cdot} \rangle$ between $\g^*$ and $\g$.
It is well known (see e.g.\ section~13.4 in \cite{marsden99itm}) that mechanical systems formulated on the cotangent bundle $T^*G$
with a left or right invariant Hamiltonian can be reduced to a system
on $\g^*$ given as
\[
     \dot{\mu} = \pm\ad^*_{\frac{\partial H}{\partial\mu}} \mu,
\]
where the negative sign is used in case of right invariance. The solution to this system preserves coadjoint orbits, which makes it natural to suggest the group action
\[
    g\cdot\mu = \Ad_{g^{-1}}^*\mu,
\]
so that the resulting Lie group integrator also respects this invariant.
For Euler's equations for the free rigid body, the Hamiltonian is left invariant and the coadjoint orbits are spheres in $\g^*\cong\R^3$.

\paragraph{Homogeneous spaces and the Stiefel and Grassmann manifolds}
The situation when $G$ acts on itself by left of right multiplication is a special case of a homogeneous space \cite{munthe-kaas97nio}, where the assumption is only that $G$ acts transitively and continuously on some manifold $M$. Homogeneous spaces are isomorphic to the quotient~$G/G_x$ where $G_x$ is the
\emph{isotropy group} for the action at an arbitrarily chosen point $x\in M$
\[
       G_x = \{ h\in G \mid h\cdot x = x\}.
\]
Note that if $x$ and $z$ are two points on $M$, then by transitivity of the action, $z=g\cdot x$ for some $g\in G$. Therefore,
whenever $h\in G_z$ it follows that $g^{-1} \cdot h \cdot g\in G_x$ so isotropy groups are isomorphic by conjugation~\cite{bryant95ait}.
Therefore the choice of $x\in M$ is not important for the construction of the quotient. For a readable introduction to this type of construction, see \cite{bryant95ait}, in particular Lecture 3.

 A much encountered example is the hypersphere
$M=S^{d-1}$ corresponding to the left action by $G=\SO(d)$,
 the Lie group of orthogonal $d\times d$ matrices with unit determinant.
 One has $S^{d-1} = \SO(d)/\SO(d-1)$. We have in fact already discussed the example of the free rigid body \eqref{eq:euler_frb} where $M=S^2$.
 
 The Stiefel manifold $\St(d,k)$ can be represented by the set of $d\times k$-matrices with orthonormal columns. An action on this set is obtained by left multiplication by $G=\SO(d)$. Lie group integrators for Stiefel manifolds are extensively studied in the literature, see e.g.\ \cite{celledoni03oti,krogstad03alc} and some applications involving Stiefel manifolds are discussed in Section~\ref{sec:dataanalysis}.
An important subclass of the homogeneous spaces is the symmetric spaces, also
obtained through a transitive action by a Lie group $G$, where $M=G/G_x$, but here one requires in addition that the isotropy subgroup is an open subgroup of the fixed point set of an involution of $G$~\cite{munthe-kaas01aos}.  
A prominent example of a symmetric space in applications is the Grassmann manifold, obtained as $\SO(d)/(\SO(k)\times \SO(d-k))$.

\paragraph{Isospectral flows} In isospectral integration one considers dynamical systems evolving on the manifold of $d\times d$-matrices sharing the same Jordan form.
Considering the case of symmetric matrices, one can use the transitive group action by $\SO(d)$
given as
\[
       g\cdot m = g m g^T.
\]
This action is transitive, since any symmetric matrix can be diagonalised by an appropriately chosen orthogonal matrix.  If the eigenvalues are distinct, then the isotropy group is discrete and consists of all matrices in $\SO(d)$ which are diagonal.

Lie group integrators for isospectral flows have been extensively studied, see for example
\cite{calvo96rkm1,calvo97nso}.  See also~\cite{celledoni01ano} for an application to the KdV equation.

\paragraph{Tangent and cotangent bundles}
For mechanical systems the natural phase space will often be the tangent bundle $TM$ as in the Lagrangian framework or the cotangent bundle $T^*M$ in the Hamiltonian framework.
The seminal paper by Lewis and Simo \cite{lewis94caf} discusses several Lie group integrators for mechanical systems on cotangent bundles, deriving methods which are symplectic, energy and momentum preserving.
Eng{\o}~\cite{engo03prk} suggested a way to generalise the Runge--Kutta--Munthe-Kaas methods into a partitioned version when $M$ is a Lie group. 
Marsden and collaborators have developed the theory of Lie group integrators from the variational viewpoint over the last two decades. See \cite{marsden01dma} for an overview.
For more recent work pertaining to Lie groups in particular, see
\cite{lee07lgv,bou-rabee09hpi,saccon09mrf}. In Section~\ref{sec:symplie} we present what we believe to be the first symplectic partitioned Lie group integrators on $T^*G$ phrased in the framework we have discussed here. Considering trivialised cotangent bundles over Lie groups is particularly attractive since there is a natural way to extend action by left multiplication from $G$ to
$G\times\g^*$ via \eqref{eq:prodGxgs}.

\subsection{Isotropy -- challenges and opportunities}
An issue which we have already mentioned a few times is that the map
$\lambda_* \colon \g \rightarrow \mathcal{X}(M)$ defined in \eqref{eq:groupaction} is not necessarily injective. This means that the choice of $f \colon M\rightarrow\g$ is not unique.
In fact, if $g \colon M\rightarrow\g$ is any map satisfying $\lambda_*(g(m))(m)=0$
for all $m\in M$, then we could replace the map $f$ by $f+g$ in \eqref{eq:genpres} without altering the vector field $F$. But such a modification of $f$ \emph{will have} an impact on the numerical schemes that we consider. This freedom in the setup of the methods makes it challenging to prove general results for Lie group methods, it might seem that some restrictions should apply to the isotropy choice for a more well defined class of schemes.
However, the freedom can of course also be taken advantage of to obtain approximations of improved quality.

\begin{figure}
\centering
 \includegraphics{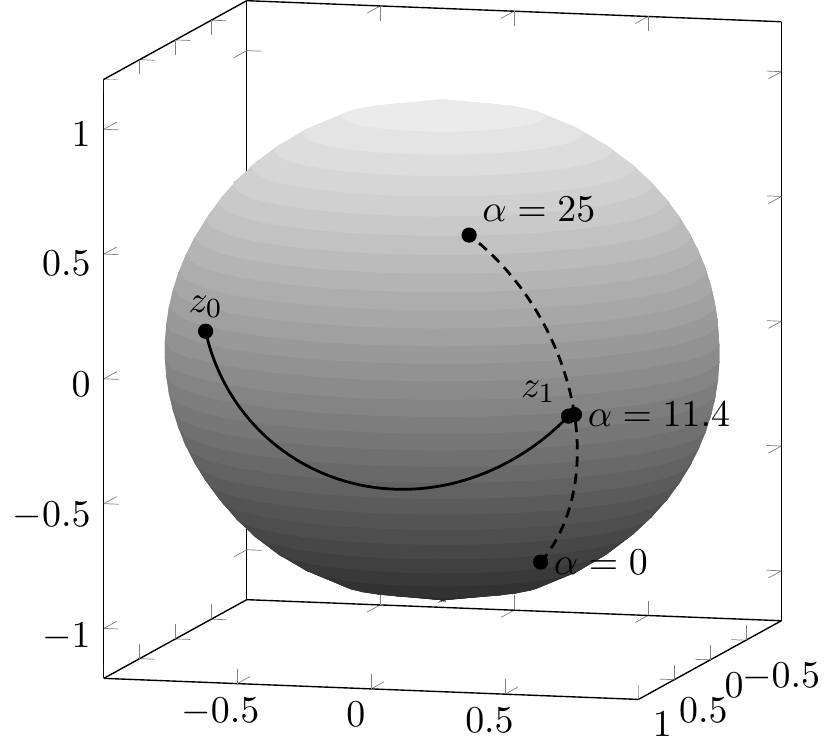}
\caption{The effect of isotropy on $S^2$ for Euler's free rigid body equations.
 The curve drawn from the initial point $z_0$ to 
$z_1$ is the exact solution, i.e.\ the momenta in body coordinates. The broken line shows the terminal points using the Lie--Euler method for $\alpha(z_0)$ (as in \eqref{eq:LieEulerS2}) varying between $0$ and $25$.
\label{fig:isotropy}}
\end{figure}

An illustrative example is the two-sphere $S^2$
acted upon linearly by the special orthogonal group $\SO(3)$. Representing elements of the 
Lie algebra $\mathfrak{so}(3)$ by vectors in $\R^3$, and points on the sphere as unit length vectors in $\R^3$, we may facilitate \eqref{eq:genpres} as
\[
     F(m)= f(m) \times m = (f(m) + \alpha(m)m)\times m,
\]
for any scalar function $\alpha \colon  S^2 \rightarrow \R$.
Using for instance the Lie--Euler method one would get
\begin{equation} \label{eq:LieEulerS2}
   m_1 = \exp(f(m_0)+\alpha(m_0)m_0) m_0,
\end{equation}
where the $\exp$ is the matrix exponential of the $3\times 3$ skew-symmetric matrix associated to a vector in $\R^3$ via the hat-map~\eqref{eq:hatmap}. Clearly the approximation depends on the choice of $\alpha(m)$. The approach of Lewis and Olver \cite{lewis08gia} was to use the isotropy to improve certain qualitative features of the solution. In particular, they studied how the orbital error could be reduced by choosing the isotropy in a clever way. In Figure~\ref{fig:isotropy} we illustrate the issue of isotropy for the 
Euler free rigid body equations. The curve drawn from the initial point~$z_0$ to 
$z_1$ is the exact solution, i.e.\ the momenta in body coordinates. The broken line shows the terminal points using the Lie--Euler method for $\alpha$ varying between $0$ and $25$.

Another potential difficulty with isotropy is the increased computational complexity when the group~$G$ has much higher dimension than the manifold $M$. This could for instance be the case with the Stiefel manifold~$\St(d,k)$ if $d\gg k$.  Linear algebra operations used in integrating differential equations on the Stiefel manifold should preferably be of complexity
$\mathcal{O}(dk^2)$.  But solving a corresponding problem in the Lie algebra $\mathfrak{so}(d)$ would typically require linear algebra operations of complexity $\mathcal{O}(d^3)$, see for example \cite{celledoni03oti} and references therein.
By taking advantage of the many degrees of freedom provided by the isotropy, it is actually possible to reduce the cost down to the required $\mathcal{O}(dk^2)$ operations as explained in for instance \cite{celledoni03lgm} and \cite{krogstad03alc}.

\section{Applications to nonlinear problems of evolution in classical mechanics}
\label{sec:appcm}

\begin{figure}
\hfil
{\includegraphics{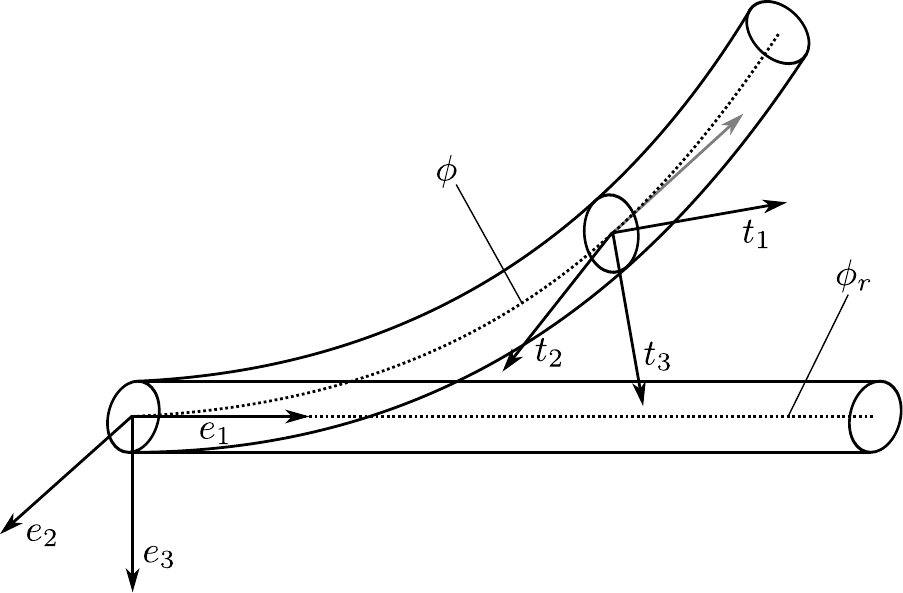}}\hfil
\caption{Geometric rod model. Here $\phi$ is the line of centroids and a cross section is identified by the frame $\Lambda=[\mathbf{t}_1,\mathbf{t}_2,\mathbf{t}_3]$, $\phi_r$ is the initial position of the line of centroids.
}
\label{fig1}
\end{figure}

The emphasis on the use of Lie groups in modelling and simulation of engineering problems in classical mechanics started in the eighties with the pioneering and fundamental work of J.C.~Simo and his collaborators.
In the case of rod dynamics, for example,  models based on partial differential equations were considered where the 
configuration  of the centreline of the rod is parametrised via arc-length, and 
the movement of a rigid frame 
attached to each of the cross sections of the rod is considered (see Figure~\ref{fig1}). 
This was first presented in a geometric context in \cite{simo85afs}. 

In robot technology, especially robot locomotion and robot grasping, the
occurrence of non-holonomically constrained models is very common.
The motion of robots equipped with wheels is not always locally
controllable, but is often globally controllable. A classical
example is the parking of a car that  cannot be moved in the
direction orthogonal to its wheels. The introduction of Lie groups
and Lie brackets to describe the dynamics of such systems, has
been considered by various authors, see for example
\cite{murray93nmp}. The design of numerical integration methods in this
context has been addressed in the paper of Crouch and Grossman,~\cite{crouch93nio}. These methods have had a fundamental impact to the successive developments in the  
field of Lie group methods. 

The need for improved understanding of non-holonomic numerical
integration has been for example advocated in
\cite{mclachlan96aso}. Recent work in this field has led to the
construction of low order non-holonomic integrators based on a
discrete Lagrange--d'Alembert principle, \cite{cortes01nhi,mclachlan06ifn}.  The use of Lie group
integrators in this context has been considered in
\cite{leok05alg,mclachlan06ifn}.

We have already mentioned the relevance of rigid body dynamics  to the numerical discretisation of rod models. There are many other research areas in which the accurate and efficient simulation of  rigid body dynamics is crucial: molecular dynamics, satellite dynamics, and celestial mechanics just to name a few, \cite{leimkuhler04shd}. 
In some of these applications, it is desirable
to produce numerical approximations which are accurate
possibly to the size of roundoff. The simulations of interest
occur over very long times and/or a large number of bodies
and this inevitably causes propagation of errors even when the
integrator is designed to be very accurate. For this reason
accurate symplectic rigid body integrators are of interest 
because they can guarantee that the roundoff error
produced by the accurate computations can stay bounded also in
long time integration.  This fact seems to be of crucial importance in 
celestial mechanics simulations, \cite{laskar04aln}. A symplectic and energy preserving Lie group integrator for the free rigid body motion was proposed in \cite{lewis94caf}. The method computes a time re-parametrisation of the exact solution.
Some recent and promising work in this
field has been  presented in \cite{mclachlan05tdm,celledoni06ets,celledoni07eco,hairer06pdm}. The control of rigid bodies with variational Lie group integrators was considered in \cite{leok05alg}.

In the next section we illustrate the use of Lie group methods in applications on a particular case study, the pipe-laying process from  ships to the bottom of the sea.

\subsection{Rigid body and rod dynamics}

\paragraph{Pipe-laying problem}

\begin{figure}
\centering
{\includegraphics{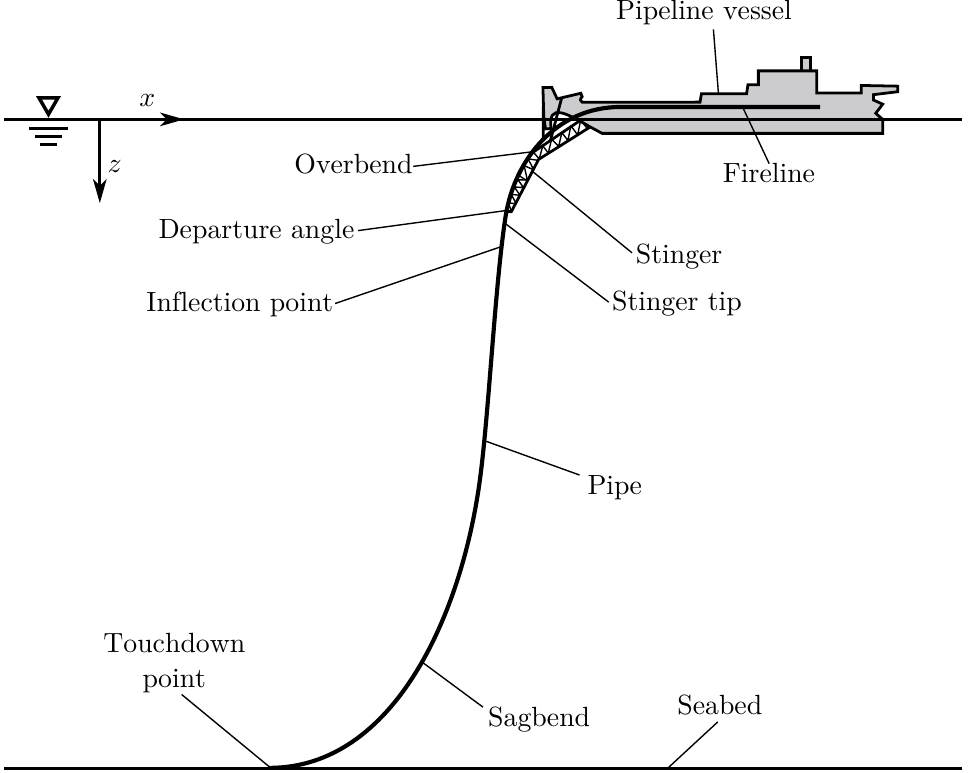}}
\caption{The pipe-laying process.
}
\label{fig2}
\end{figure}

The simulation of deep-water risers, pipelines and drill rigs
requires the use of models of long and thin beams subject to large
external forces. These are 
complicated nonlinear systems with highly oscillatory
components. 
We are particularly interested in the correct and
accurate simulation of the pipe-laying process
 from
ships on the bottom of the sea, see Figure~\ref{fig2}. The problem comprises the modelling of two interacting structures: a long and thin pipe (modelled as a rod) and a vessel (modelled as a rigid body). The system is subject to environmental forces (such as sea and wind effects). The control parameters for this problem are the vessel position and velocity, the pay-out speed and the pipe tension while the control objectives consist in determining the touchdown position of the pipe as well as ensuring the integrity of the pipe and to avoid critical deformations, \cite{jensen10anp,safstrom09mas}.

The vessel rigid body equations determine the boundary conditions of the rod. They
are expressed in six degrees of freedom
as
\begin{equation*}\label{eq:vesseleq}
M\dot{\boldsymbol\nu} + C(\boldsymbol\nu)\boldsymbol\nu
 + D(\boldsymbol\nu)\boldsymbol\nu + g(\boldsymbol\eta) = 
\boldsymbol\tau,
\end{equation*}
where $M$ is the system inertia matrix, $C(\boldsymbol\nu)$ the Coriolis-centripetal
matrix, $D(\boldsymbol\nu)$ the damping matrix, $g(\boldsymbol\eta)$ the vector of gravitational
and buoyancy forces and moments, and $\boldsymbol\tau$ the vector of control inputs and
environmental disturbances such as wind, waves and currents (see \cite{perez07kmf} for details). The vector $\boldsymbol\nu$ contains linear and angular velocity and $\boldsymbol\eta$ is the position vector.
It has been shown in \cite{jensen10anp} that the rigid body vessel equations are input-output passive. 

The equations can be integrated numerically with a splitting and composition technique where the vessel equations are split into a free rigid body part and a damping and control part. The free rigid body equations can be solved with a method proposed in \cite{celledoni06ets} where the angular momentum is accurately and efficiently computed by using Jacobi elliptic functions, the attitude rotation is obtained using a Runge--Kutta--Munthe-Kaas Lie group method, and the control and damping part is solved exactly.

Simulations of the whole pipe-lay problem with local parametrisations of the pipe and the vessel based on Euler angles have been obtained in \cite{jensen10anp}.

\paragraph{Rod dynamics}
At fixed time 
each cross section of the pipe is the result of a rigid rotation in space of a reference cross section, and analogously, for each fixed value of the space variable the corresponding cross section evolves in time as a forced rigid body, see Figure~\ref{fig1}.
In absence of external forces the equations are
\begin{align*}
\rho_A \partial_{tt} \phi &= \partial_S\mathbf{n},\qquad S\in[0,L],\quad t\ge0,\\
\partial_t \pi+ (I_{\rho}^{-1} \pi)\times \pi &= \partial_S\mathbf{m}+(\partial_S \phi)\times \mathbf{n},
\end{align*}
here $\phi=\phi(S,t)$ is the line of centroids of the rod,  $\mathbf{m}$ and $\mathbf{n}$ are the stress couple and stress resultant, $\mathbf{\pi}$ is the angular momentum density, $I_{\rho}$ is the inertia in spatial coordinates, and $\rho_A=\rho_A(S)$ is the mass per unit length of the rod (see \cite{simo87otd} and \cite{celledoni10aha}).
The kinematic equations for the attitude rotation matrix are
\[
	\partial_t\Lambda=\hat{\mathbf{w}} \Lambda, \qquad \partial_S\Lambda=\hat{\mathbf{M}} \Lambda,
\]
 where $\Lambda(S,t)=[\mathbf{t}_1,\mathbf{t}_2,\mathbf{t}_3]$, $I_{\rho}^{-1}\pi=\mathbf{w}$, $\mathbf{M}=C_2 \Lambda^T\mathbf{m}$ and $C_2$ is a constant diagonal matrix. 
 We denote by ``$\hat{\quad}$"  the hat-map identifying $\R^3$ with $\mathfrak{so}(3)$:
\begin{equation} \label{eq:hatmap}
    \mathbf{v}=\begin{bmatrix}
            v_1\\ v_2\\ v_3
        \end{bmatrix} \mapsto 
        \hat{\mathbf{v}}= \begin{bmatrix*}[r]
            0 &-v_3&v_2\\v_3&0&-v_1\\ -v_2&v_1&0
        \end{bmatrix*}.
\end{equation}
With no external forces one assumes pure displacement boundary conditions providing $\phi$ and $\Lambda$ on the boundaries $S=0$ and $S=L$.

In \cite{simo87otd}, partitioned Newmark
integrators, of Lie group type, and of moderate order were considered for this problem. 
While classical Newmark methods  are variational integrators and as such are symplectic when applied to Hamiltonian systems \cite{marsden01dma}, the particular approach of \cite{simo87otd} involves the use of exponentials for the parametrisation of the Lie group $\SO(3)$, and the geometric properties of this approach are not trivial to analyse. 
Moreover, since the model is a partial differential equation, space and time discretisations should be designed so that the overall discrete equations admit solutions and are stable. It turns out that conventional
methods perform poorly on such problem in long time simulations. 
To obtain stable methods reliable in long-time simulation, an energy-momentum method was proposed for the rod problem in \cite{simo95ndo}. 
Later, this line of thought has been further developed in \cite{romero01aof}.
The Hamiltonian formulation of this model allows one to derive natural structure preserving discretisations into systems of coupled rigid bodies \cite{leyendecker06oem}.

Following the geometric space-time integration procedure proposed in \cite{frank04gst}, a multi-Hamiltonian formulation\footnote{For a definition of the multi-symplectic structure of Hamiltonian partial differential equations, see \cite{bridges01msi}.} of these equations has been proposed in \cite{celledoni10aha}, using the Lie group of Euler parameters. The design of corresponding multi-symplectic Lie group discretisations is still under investigation. 


\section{Applications to problems of data analysis and statistical signal processing}
\label{sec:dataanalysis}


The solution of the many-body Schr{\"o}dinger eigenvalue problem
\begin{equation}
\label{eq:schr}
\hat{\mathbf{H}}\Psi =E\Psi,
\end{equation}
where the so called electronic ground state (the smallest eigenstate) is sought, is an important problem of computational chemistry. The main difficulty is the curse of dimensionality. Since $\hat{\mathbf{H}}$ is a differential operator in several space dimensions, a realistic simulation of (\ref{eq:schr}) would require the numerical discretisation and solution of a partial differential equation in several space dimensions. The number of space dimensions  grows with the number of electrons included in the simulation.

The eigenvalue problem admits an alternative variational formulation. Instead of looking for the smallest eigenvalue and eigenfunction of the Schr{\"o}dinger equation, one minimises directly the ground state energy. 

After appropriate spatial discretisation, 
the problem becomes a minimisation problem on a Riemannian manifold
$\mathcal{M}$, 
\begin{equation}\label{genopt}
{\min_{x\in\mathcal{M}}}\,
\phi(x),
\end{equation}
where
$\phi \colon \mathcal{M}\rightarrow \R$  is a smooth discrete energy function to be minimised on $\mathcal{M}$. The discrete energy $\phi$ considered  here is the so called  Kohn--Sham energy, \cite{baarman12dma}. For a related application of Lie group techniques in quantum control, see \cite{degani09qcw}.

The general optimisation problem giving rise to (\ref{genopt}) appears in several
applied fields, ranging from engineering to applied physics and
medicine. 
Some specific examples are principal
component/subspace analysis, eigenvalue and generalised eigenvalue
problems, optimal linear compression, noise reduction, signal
representation and blind source separation.  

\subsection{Gradient-based optimisation on Riemannian manifolds}

Let $\mathcal{M}$ be a Riemannian  manifold with metric $\langle\cdot, \cdot
\rangle_{\mathcal{M}}$ and $\phi \colon \mathcal{M}\rightarrow \R$  be a smooth cost function to be minimised on $\mathcal{M}$.
We want to solve (\ref{genopt}). 
The
optimisation method based on gradient flow -- written for the
minimisation problem only, for the sake of easy reading -- consists in
setting up the differential equation on the manifold,
\begin{equation}\label{eqdiff}
\dot{x}(t)=-\grad \phi\big(x(t)\big),
\end{equation}
with appropriate initial condition $x(0)=x_0\in\mathcal{M}$. The equilibria
of equation (\ref{eqdiff}) are the critical points of the function
$\phi$. 
In the above equation, the symbol $\grad \phi$
denotes the Riemannian gradient of the function $\phi$ with respect
to the chosen metric. Namely, $\grad \phi(x)\in T_x\mathcal{M}$ and
$T_x \phi(v)=\langle \grad \phi (x),v\rangle_{\mathcal{M}}$
for all $v\in T_x\mathcal{M}$. 

The solution of (\ref{eqdiff}) on $\mathcal{M}$ may be locally  expressed in terms of a
curve on the tangent space $T_{x_0}\mathcal{M}$ using a
retraction map $\mathcal{R}$. 
Retractions are tangent space parametrisations of $\mathcal{M}$, and allow us to write
\[x(t)=\mathcal{R}_{x_0}(\sigma(t)), \quad \sigma(t)\in T_{x_0}\mathcal{M}, \quad t\in [0,t_f],\]
for small enough $t_f$, see  \cite{shub86sro} for a precise definition.

In most applications of interest, see for example \cite{brockett91dst,helmke94oad}, $\mathcal{M}$ is a matrix manifold endowed with a Lie group action and there is a natural way to define a metric and a retraction.
In fact, let $\mathcal{M}$ be a manifold acted upon 
by a Lie group $G$, with a locally transitive group action
$\Lambda(g,x)=\Lambda_x(g)$. Let us also consider a coordinate map $\psi$,
\[
\psi \colon \g \rightarrow G, \quad \text{and} \quad  \rho_x \mdef T_0(\Lambda_x\circ
\psi).
\]

One can prove that if there exists a linear map
$a_x \colon T_x\mathcal{M}\rightarrow \g$ such that $\rho_x\circ
a_x =\Id_{T_x\mathcal{M}}$, then $\mathcal{R}_x$, given by
\[ 
\mathcal{R}_x(v) \mdef 
(\Lambda_x \circ \psi
\circ a_x )(v),
\]
is a retraction, see \cite{celledoni03oti}. The existence of $a_x$ is guaranteed, at least locally, by the transitivity of the action and the fact that $\psi$ is a local diffeomorphism. The approach is analogous to the one described for differential equations in section~\ref{GeneralizingRK}. Therefore, we can construct
retractions using any coordinate map from the Lie algebra $\g$ to
the group. 

Any metric on $\mathfrak{g}$, $\langle\cdot ,\cdot\rangle_{\mathfrak{g}}$ induces a metric on $\mathcal{M}$ by
\[\langle v_x,w_x\rangle_{\mathcal{M}}=\langle a_x(v_x) ,a_x(w_x)\rangle_{\mathfrak{g}}.\]
Also, we may define the image of the tangent space under the map
$a_x$:
\[
\mathfrak{m}_x \mdef a_x(T_x\mathcal{M})\subset \g.
\]
The set $\ensuremath{\mathfrak{m}}_x$ is a linear subspace of the
Lie algebra $\g$, often of lower dimension. Parametrisations of the solution of (\ref{eqdiff}) involving the whole Lie algebra are in general more computationally intensive than those restricted to $\ensuremath{\mathfrak{m}}_x$, but, if the isotropy is chosen suitably, they might lead to methods which converge faster to the optimum.

For the sake of illustration, we
consider the minimisation on a two-dimensional torus $T^2=S^1\times
S^1$. 
Here we denote by $S^1$ the circle, i.e.\
\begin{gather*}
S^1=\{g(\alpha) \mathbf{e}_1 \in \R^2 \mid g(\alpha)\in \SO(2)\}, \\
g(\alpha)=\exp(\alpha E), \quad E=\begin{bmatrix*}[r]
                                                                                               0 & -1\\
                                                                                              1 & 0\\
\end{bmatrix*}, \quad 0 \le \alpha < 2\pi,
\end{gather*}
where
$\mathbf{e}_1$ is the first canonical vector and $\SO(2)$ is the 
commutative Lie group of planar rotations.
Any element in $T^2$ is of the form 
\[x_0\in T^2,\quad x_0=(g(\theta)\mathbf{e}_1,g(\varphi)\mathbf{e}_1), 
\quad g(\theta), g(\varphi)\in \SO(2).\]
The Lie group acting on $T^2$ is $\SO(2)\times \SO(2)$, its corresponding Lie algebra is 
$\mathfrak{so}(2)\times \mathfrak{so}(2)$, which has dimension $d=2$ and basis $\{ (E,O), (O,E)\}$, where $O$ is the 
zero element in $\mathfrak{so}(2)$.  

The Lie group action is
\[\Lambda_{x_0}(h_1,h_2)=(h_1 g(\theta)\mathbf{e}_1,h_2 g(\varphi)\mathbf{e}_1),\qquad (h_1,h_2)\in \SO(2)\times \SO(2),\]
and $\psi=\exp$.  Any $v_{x_0}\in T_{x_0}T^2$ can be written as 
\[v_{x_0}=(\alpha E\mathbf{e}_1,\beta E\mathbf{e}_1),\]
for some $\alpha, \beta \in\R,$
so
\[a_{x_0}(v_{x_0})=(\alpha E,\beta E).\]

\begin{figure}
\centering
\ifdraft{\tikz \draw (0cm,0cm) rectangle (12cm,10cm);}{\input{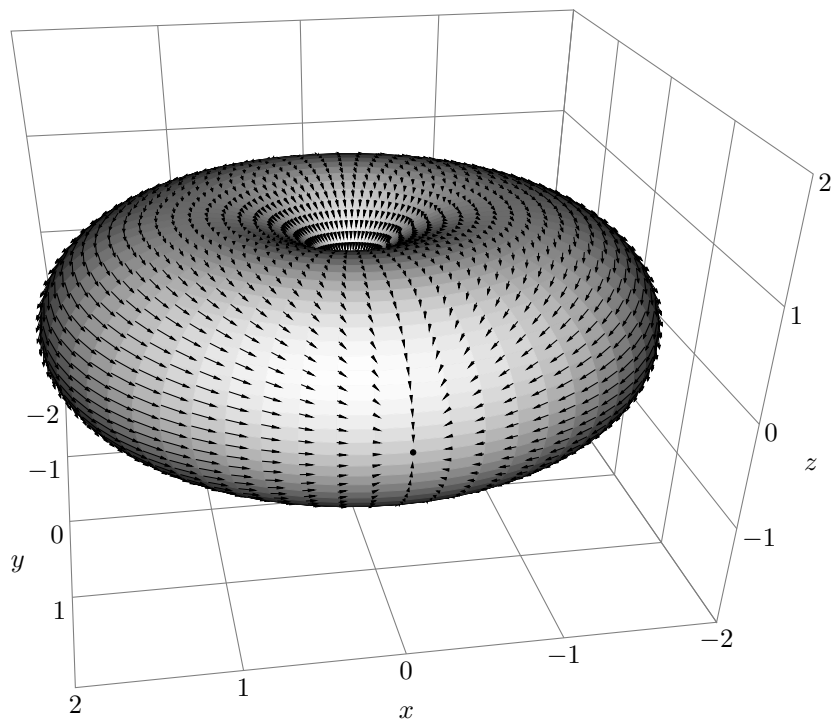}}
\caption{The gradient vector field of the cost function
$\phi (g(\theta)\mathbf{e}_1,g(\varphi)\mathbf{e}_1)= ((1+\mathbf{e}_1^Tg(\theta)\mathbf{e}_1) \mathbf{e}_2^Tg(\varphi)\mathbf{e}_1-8)^2$ on the torus.
The vector field points towards the two minima, the global minimum is marked with a black spot in the middle of the picture.
}
\label{fig0}
\end{figure}

Assume the cost function we want to minimise is simply the distance from a
fixed plane in $\R^3$, say the plane with equation $y=8$. This gives
\begin{equation*}
\label{cost1}
\phi(g(\theta)\mathbf{e}_1,g(\varphi)\mathbf{e}_1)=\lvert(1+\cos (\theta)) \sin (\varphi)-8\rvert,
\end{equation*} 
and the
minimum is attained in $\theta =0$ and $\varphi=\pi/2$.\footnote{
We have used a parameterisation of $T^2$ in $\R^3$ in angular coordinates,  obtained applying the following mapping
\[
(g(\theta)\mathbf{e}_1,g(\varphi)\mathbf{e}_1) \mapsto \left\{
\begin{alignedat}{2}
   x &= (1+\mathbf{e}_1^Tg(\theta)\mathbf{e}_1) \mathbf{e}_1^Tg(\varphi)\mathbf{e}_1 & &= (1+\cos (\theta)) \cos (\varphi),\\
   y &= (1+\mathbf{e}_1^Tg(\theta)\mathbf{e}_1) \mathbf{e}_2^Tg(\varphi)\mathbf{e}_1 & &= (1+\cos (\theta)) \sin (\varphi),\\
   z &= \mathbf{e}_2^Tg(\theta)\mathbf{e}_1 & &= \sin (\theta),
\end{alignedat} \right.
\]
with $0 \le \theta , \varphi < 2\pi$.
This is equivalent to the composition of two planar rotations and one
translation in $\R^3.$}
In Figure~\ref{fig0} we plot $-\grad \phi$, 
the negative gradient vector field for the given cost function. The Riemannian metric we used is
\[
\langle(\alpha_1E\mathbf{e}_1,\beta_1E\mathbf{e}_1),(\alpha_2E\mathbf{e}_1,\beta_2E\mathbf{e}_1)\rangle_{T^2}=\alpha_1\alpha_2+\beta_1\beta_2,
\]
and $(\alpha_1E\mathbf{e}_1,\beta_1E\mathbf{e}_1)\in T_{(\mathbf{e}_1,\mathbf{e}_1)}T^2$. This metric can be easily interpreted as a metric on the Lie algebra~$\mathfrak{g}=\mathfrak{so}(2)\times \mathfrak{so}(2)$:
\[
\langle(\alpha_1E,\beta_1E),(\alpha_2E,\beta_2E)\rangle_{\mathfrak{g}}=\alpha_1\alpha_2+\beta_1\beta_2.\]
At the point 
$p_0=( g(\theta_0)\mathbf{e}_1,  g(\varphi_0)\mathbf{e}_1)\in T^2$, 
the gradient vector field can be represented by
\[
(\gamma E g(\theta_0)\mathbf{e}_1,\delta E g(\varphi_0)\mathbf{e}_1),
\]
where $\gamma$ and $\delta$ are real values given by
\[
\gamma=-C\sin(\theta_0)\sin(\varphi_0),\qquad \delta=C(1+\cos(\theta_0))\cos(\varphi_0),\]
and
\[C=2((1+\cos(\theta_0)\sin(\varphi_0)-8).
\]

Gradient flows are not the only type of differential equations which can be used to solve optimisation problems on manifolds. Alternative equations have been proposed in the context of neural networks \cite{celledoni04nlb, celledoni08dmf}. Often they arise naturally as the Euler--Lagrange equations of a variational problem.

\subsection{Principal component analysis}

Data reduction techniques are statistical signal processing methods
that aim at providing efficient representations of data. A
well-known data compression technique consists of mapping a
high-dimensional data space into a lower dimensional representation
space by means of a linear transformation. It requires the
computation of the data covariance matrix and then the application
of a numerical procedure to extract its eigenvalues and the
corresponding eigenvectors. Compression is then obtained by representing the signal in a basis consisting only 
of
those 
eigenvectors associated with the most significant
eigenvalues. 

In particular, principal component analysis (PCA) is a second-order
adaptive statistical data processing technique that helps removing
the second-order correlation among given random signals. Let us
consider a stationary multivariate random process $x(t)\in\R^{n}$
and suppose its covariance matrix $A=E[(x-E[x])(x-E[x])]^{T}]$
exists and is bounded. Here the symbol $E[\cdot]$ denotes statistical expectation.
If $A\in\R^{n\times n}$ is not diagonal, then the
components of $x(t)$ are statistically correlated. One can remove this 
redundancy by partially diagonalising $A$, i.e.\ computing 
the operator $F$ formed by the eigenvectors of the matrix $A$
corresponding to its largest eigenvalues. This is possible since the covariance matrix $A$
is symmetric (semi) positive-definite, and $F\in\St(n,p)$. 

To compute $F$ and the corresponding $p$ eigenvalues of the $n\times n$ symmetric and positive-definite matrix~$A$, we consider the maximisation of the function
\[
	\phi(X)=\frac{1}{2}\operatorname{trace}(X^TAX),
\]
on the Stiefel manifold, and solve numerically the corresponding gradient flow with a Lie group integrator.

As a consequence
the new random signal
defined by $y(t) \mdef F^{T}(x(t)-E[x(t)])\in\R^{p}$ has uncorrelated
components, with $p\leq n$ properly selected. 
The
component signals of $y(t)$ are the so called \emph{principal components
of the signal} $x(t)$, and their relevance is proportional to the corresponding
eigenvalues $\sigma_{i}^2=E[y_{i}^2]$ which here are
arranged in descending order ($\sigma_{i}^{2}\geq\sigma_{i+1}^{2}$).

Thus, the data stream $y(t)$ is a compressed version of the
data stream $x(t)$. After the reduced-size data has been processed
(i.e.\ stored, transmitted), it needs to be recovered, that is,
it needs to be brought back to the original structure. However, the
principal-component-based data reduction technique is not lossless,
thus only an approximation $\hat{x}(t)\in\R^{n}$ of the original
data stream may be recovered. An approximation of $x(t)$ is given by $\hat{x}(t)=Fy(t)+E[x]$.
Such approximate data stream minimises the reconstruction error~$E[\lVert x-\hat{x}\rVert_2^{2}]=\sum_{i=n+1}^{p}\sigma_{i}^2$.

For a scalar or a vector-valued random variable $x\in\R^n$ endowed with a probability
density function~$p_x \colon x\in\R^n\rightarrow p_x(x)\in\R$, the expectation of a function
$\beta \colon \R^n\rightarrow \R$ is defined as
\[
E[\beta] \mdef \int_{\R^n}\beta(x)p_x(x) \, \d^nx.
\]
Under the hypothesis that the signals whose expectation is to be
computed are ergodic, the actual expectation (ensemble average) may
be replaced by temporal-average on the basis of the available
signal samples, namely
\[
E[\beta]\approx \frac{1}{T}\sum_{t=1}^T\beta(x(t)).
\]

\subsection{Independent component analysis}

An interesting example of a problem that can be tackled via
statistical signal processing is the \emph{cocktail-party problem}.
Let us suppose $n$ signals $x_1(t),\dots ,x_n(t)$ were recorded from
$n$ different positions in a room where there are $p$ sources or speakers.
Each recorded signal is a linear mixture of the voices of the sources~$s_1(t),\dots , s_p(t)$, namely
\begin{align*}
x_1(t) &= a_{1,1} s_1(t) +\dotsb +a_{1,p} s_p(t),\\
&\vdotswithin{=} \\
x_n(t) &= a_{n,1} s_1(t) +\dotsb +a_{n,p} s_p(t),
\end{align*}
where the $n p$ coefficients $a_{i,j}\in\R$ denote the mixing
proportions. The mixing matrix $A=(a_{i,j})$ is unknown. The cocktail party problem consists in estimating
signals $s_1(t),\dots , s_p(t)$ from only the knowledge of their
mixtures $x_1(t), \dots ,x_n(t)$. The main assumption on the source
signals is that $s_1(t), \dots , s_p(t)$ are statistically independent.
This problem can be solved using independent component analysis (ICA).

Typically, one has $n>p$, namely, the number of observations exceeds the number of actual sources.
Also, a typical assumption is that the source signals are spatially white, which means $E[ss^T]=I_p$, the $p \times p$~identity matrix.
The aim of independent component analysis is to find estimates $y(t)$ of signals in $s(t)$ by
constructing a de-mixing matrix $W\in\R^{n\times p}$ and by computing $y(t) \mdef W^Tx(t)$.
Using statistical signal processing methods, 
the problem is reformulated into an optimisation problem on homogeneous manifolds for finding the de-mixing
matrix $W$. 

The geometrical structure of the parameter space in ICA comes from a signal
pre-processing step named \emph{signal whitening}, which is operated on the
observable signal $x(t)\rightarrow \tilde{x}(t)\in\R^p$ in such a way that the
components of the signal $\tilde{x}(t)$ are uncorrelated and have variances equal to $1$,
namely $E[\tilde{x}\tilde{x}^T]=I_p$.
 This also means that redundant
observations are eliminated and the ICA problem is brought back to the smallest dimension $p$. This can be done by computing
$E[xx^T]=VDV^T$,
with $V\in\St(n,p)$ and $D\in\R^{p\times p}$ diagonal invertible. Then
%
\[
\tilde{x}(t) \mdef D^{-\frac{1}{2}}V^Tx(t),
\]
and with $\tilde{A} \mdef D^{-\frac{1}{2}}V^TA$ we have
$E[\tilde{x}\tilde{x}^T]=\tilde{A}E[ss^T]\tilde{A}^T=\tilde{A}\tilde{A}^T=I_p$. 

After observable signal
whitening, the de-mixing matrix may be searched for such that it solves the optimisation problem
\begin{equation*}\label{pri1}
{\max_{W\in \LieO(p)}}\,\phi(W).
\end{equation*}
 As explained, after whitening, the number of projected observations in the signal
$\tilde{x}(t)$ equals the number of sources. However, in some applications it is known
that not all the source signals are useful, so it is sensible to analyse only
a few of them. 
In these cases, if we denote by $\overline{p}\ll p$ the actual number of independent
components that are sought after, the appropriate way to cast the optimisation problem for
ICA is
\begin{equation*}\label{pri2}
\max_{W\in \St(n,\overline{p})} \phi(W), \quad \text{with} \quad \overline{p}\ll p.
\end{equation*}
The corresponding gradient flows obtained in this case are differential equations on the orthogonal group or on the Stiefel manifold, and can be solved numerically by Lie group integrators.

As a possible principle for reconstruction, the maximisation or minimisation of
non-Gaussianity is viable. It is based on the notion that the sum
of independent random variables has distribution closer
to Gaussian than the distributions of the original random variables.
A measure of non-Gaussianity is the kurtosis, defined for a scalar signal $z\in\R$ as
\[
\operatorname{kurt}(z) \mdef E[z^4]-3E^2[z^2].
\]
If the random signal $z$ has unitary variance, then the kurtosis computes as
$\operatorname{kurt}(z)=E[z^4]-3$. Maximising or minimising kurtosis is thus a possible way
of estimating independent components from their linear mixtures, see \cite{celledoni08dmf} and references therein for more details. 

\subsection{Computation of Lyapunov exponents}

The Lyapunov exponents of a continuous dynamical system
$\dot x=F(x)$, $x(t)\in \R^n$, provide a
qualitative measure of its complexity. They are numbers related to the linearisation $A(t)$ of $\dot x=F(x)$ along a trajectory~$x(t)$.
Consider the solution $U$ of the matrix problem
\[
\dot{U}=A(t)U, \qquad U(0)=U_0, \qquad U(t) \in \R^{n\times n}.
\]
The logarithms of the eigenvalues of the matrix
\[
\Lambda=\lim_{t\rightarrow \infty}\left(U(t)^{T}U(t)\right)^{\frac{1}{2t}},
\]
are the Lyapunov exponents for the given dynamical system. In
\cite{dieci95coa} a procedure for computing just
$k$ of the $n$ Lyapunov exponents of a dynamical system is presented. The exponents are computed by solving an
initial value problem on $\St (n,k)$ and
computing a quadrature of the diagonal entries of a $k\times k$ matrix-valued function.  The
initial value problem is defined as follows:
\[
\dot{Q}=(A-QQ^{T}A+QSQ^{T})Q,
\]
with random initial value in $\St(n,k)$ and 
\[
S_{k,j}= \begin{cases}
                (Q^TAQ)_{k,j}, & k>j,\\
                  0,           & k=j,\\
               -(Q^TAQ)_{j,k}, & k<j,
              \end{cases}
              \qquad k,j=1,\dots ,p.
\]

It can be shown that the $i$-th Lyapunov exponent $\lambda_i$ can be
obtained as
\begin{equation}
\label{int}
\lambda_i=\lim_{t\rightarrow
    \infty}\frac{1}{t}\int_0^tB_{i,i}(s) \, \d s, \quad i=1,\dots ,k,
\end{equation}
and
\[
B=Q^TAQ-S.
\]
 One could use for example the trapezoidal rule to
approximate the integral (\ref{int}) and compute $\lambda_i$ ($i=1,\dots
,k$). We refer to \cite{dieci95coa} for further
details on the method, and to \cite{celledoni02aco} for the use of Lie group integrators on this problem. Lie group methods for ODEs on Stiefel manifolds have also been considered in \cite{celledoni03oti,krogstad03alc,celledoni03cfl}.

We have here presented a selection of applications that can be dealt with by solving differential equations on Lie groups and homogeneous manifolds. For these problems, Lie group integrators are a natural choice. We gave particular emphasis to problems of evolution in classical mechanics and problems of signal processing. This is by no means an exhaustive survey; other interesting areas of application are for example problems in vision and medical imaging, see for instance  \cite{xu12aol,lee07gds}.

\section{Symplectic integrators on the cotangent bundle of a Lie group}
\label{sec:symplie}
In this section we shall assume that the manifold is the cotangent bundle $T^*G$ of a 
Lie group $G$. Let $R_g\colon G\rightarrow G$ be the right multiplication operator such that
$R_g(h)=h \cdot g$ for any $h\in G$. The tangent map of $R_g$ is denoted $R_{g*} \mdef TR_g$.
Any cotangent vector $p_g\in T_g^*G$ can be associated to $\mu\in\g^*$ by right trivialisation as follows: Write $v_g\in T_gG$ in the form $v_g = R_{g*}\xi$ where $\xi\in\g$, so that
$\langle p_g, v_g\rangle=\langle p_g, R_{g*}\xi\rangle = \langle R_g^* p_g, \xi\rangle$, where we have used $R_g^*$ for the dual map of $R_{g*}$, and $\langle{\cdot},{\cdot}\rangle$ is a duality pairing.
We therefore represent $p_g\in T_g^*G$ by $\mu=R_g^* p_g\in\g^*$. Thus, we may use as phase space $G\times\g^*$ rather than $T^*G$. For applying Lie group integrators we need a transitive group action on $G\times\g^*$ and this can be achieved by lifting the group structure of $G$ and using left multiplication in the extended group. The semidirect product structure on
$\mathbf{G} \mdef G\ltimes\g^*$ is defined as
\begin{equation} \label{eq:prodGxgs}
      (g_1, \mu_1)\cdot (g_2, \mu_2) = (g_1\cdot g_2, \mu_1 + \Ad_{g_1^{-1}}^*\mu_2),
\end{equation}
where the coadjoint action $\Ad^{*}$ is defined in \eqref{eq:coadjoint-action}.
Similarly, the tangent map of right multiplication extends as
\[
TR_{(g,\mu)}(R_{h*}\,\zeta, \nu) = (R_{hg*}\;\zeta, \nu-\ad_\zeta^*\,\Ad_{h^{-1}}^*\mu),\quad
g,h \in G,\ \zeta\in\g,\ \mu,\nu\in\g^*.
\]
Of particular interest is the restriction of $TR_{(g,\mu)}$ to $T_e\mathbf{G}\cong \g \times \g^*$,
\[
    T_eR_{(g,\mu)}(\zeta,\nu) = (R_{g*}\zeta, \nu - \ad_\zeta^*\mu). \]
The natural symplectic form on
$T^*G$ (which is a differential two-form) is defined as
\[
\Omega_{(g,p_g)}((\delta v_1, \delta \pi_1),(\delta v_2,\delta\pi_2))
=\langle \delta \pi_2, \delta v_1\rangle - \langle \delta \pi_1, \delta v_2\rangle,
\]
and by right trivialisation it may be pulled back to $\mathbf{G}$ and then takes the form
\begin{equation} \label{eq:sympform}
      \omega_{(g,\mu)}( (R_{g*} \xi_1, \delta\nu_1), (R_{g*}\xi_2, \delta\nu_2))
      = \langle\delta\nu_2,\xi_1 \rangle - \langle\delta\nu_1, \xi_2 \rangle - \langle\mu, [\xi_1,\xi_2]\rangle, \qquad \xi_1,\xi_2 \in \g.
\end{equation}
The presentation of differential equations on $T^*G$ as in \eqref{eq:genpres} is now achieved via the action by left multiplication, meaning that any vector field $F\in\mathcal{X}(\mathbf{G})$ is expressed by means of a map $f\colon \mathbf{G} \rightarrow T_e \mathbf{G}$,
\begin{equation} \label{eq:Fpres}
        F(g,\mu) = T_e R_{(g,\mu)} f(g,\mu) = (R_{g*} f_1, f_2-\ad_{f_1}^*\mu),
\end{equation}
where $f_1=f_1(g,\mu)\in\g$, $f_2=f_2(g,\mu)\in\g^*$ are the two components of $f$.
We are particularly interested in the case that $F$ is a Hamiltonian vector field
which means that $F$ satisfies the relation
\begin{equation} \label{eq:iF}
     \mathrm{i}_{F}\omega = \d H,
\end{equation}
for some Hamiltonian function $H\colon T^*G\rightarrow\R$ and $\mathrm{i}_F$ is the interior product defined as $\mathrm{i}_F \omega (X) \mdef \omega(F, X)$ for any vector field $X$.
From now on we let $H \colon \mathbf{G} \to \R$ denote the trivialised Hamiltonian.
A simple calculation using 
\eqref{eq:sympform}, \eqref{eq:Fpres} and \eqref{eq:iF} shows that the corresponding map $f$ for such a Hamiltonian vector field is
\[
     f(g,\mu) = \left(\frac{\partial H}{\partial\mu}(g,\mu), -R_g^*\frac{\partial H}{\partial g}(g,\mu)\right).
\]
We have come up with the following family of symplectic Lie group integrators on
$G\times\g^*$
\begin{align*}
(\xi_i, \bar{n}_i) &= h f(G_i, M_i),\qquad n_i  = \Ad_{\exp(X_i)}^* \bar{n}_i,\quad i=1,\ldots s, \\
(g_1, \mu_1) &= \exp\Bigl(Y, (\dexp_{Y}^{-1})^{*}\sum_{i=1}^s b_i n_i\Bigr)\cdot (g_0, \mu_0), \\
Y&= \sum_{i=1}^s b_i \xi_i,\qquad  X_i = \sum_{j=1}^s a_{ij} \xi_j,\quad i=1,\ldots,s, \\
G_i &= \exp(X_i)\cdot g_0,\quad i=1,\ldots, s, \\
M_i &= \dexp_{-Y}^* \mu_0 + \sum_{j=1}^s \left(b_j\dexp_{-Y}^* - \frac{b_ja_{ji}}{b_i}\dexp_{-X_j}^*\right) n_j,\quad i=1,\ldots,s.
\end{align*}
Here, $a_{ij}$ and $b_i$ are coefficients where it is assumed that $\sum_{i=1}^s b_i = 1$ and that $b_i\neq 0,\ 1\leq i\leq s$. The symplecticity of these schemes is a consequence of their derivation from a variational principle, following ideas similar to that of \cite{bou-rabee09hpi} and \cite{saccon09mrf}.  
One should be aware that order barriers for this type of schemes may apply, and that further stage corrections may be necessary to obtain high order methods.

\paragraph{Example, the $\theta$-method for a heavy top}
Let us choose $s=1$ with coefficients $b_1=1$ and $a_{11}=\theta$, i.e.\ the RK coefficients of the well known $\theta$-method.
Inserting this into our method and simplifying gives us the method
\begin{align*}
	(\xi, \bar n) &= h f\bigl(\exp(\theta \xi) \cdot g_0, \dexp_{-\xi}^{*} \mu_0 + (1 - \theta) \dexp_{-(1-\theta)\xi}^{*} \bar n\bigr), \\
	(g_1, \mu_1) &= (\exp(\xi), \Ad_{\exp(-(1 -\theta) \xi)}^{*} \bar n) \cdot (g_0, \mu_0).
\end{align*}
\begin{figure}\centering
\includegraphics[width=\textwidth]{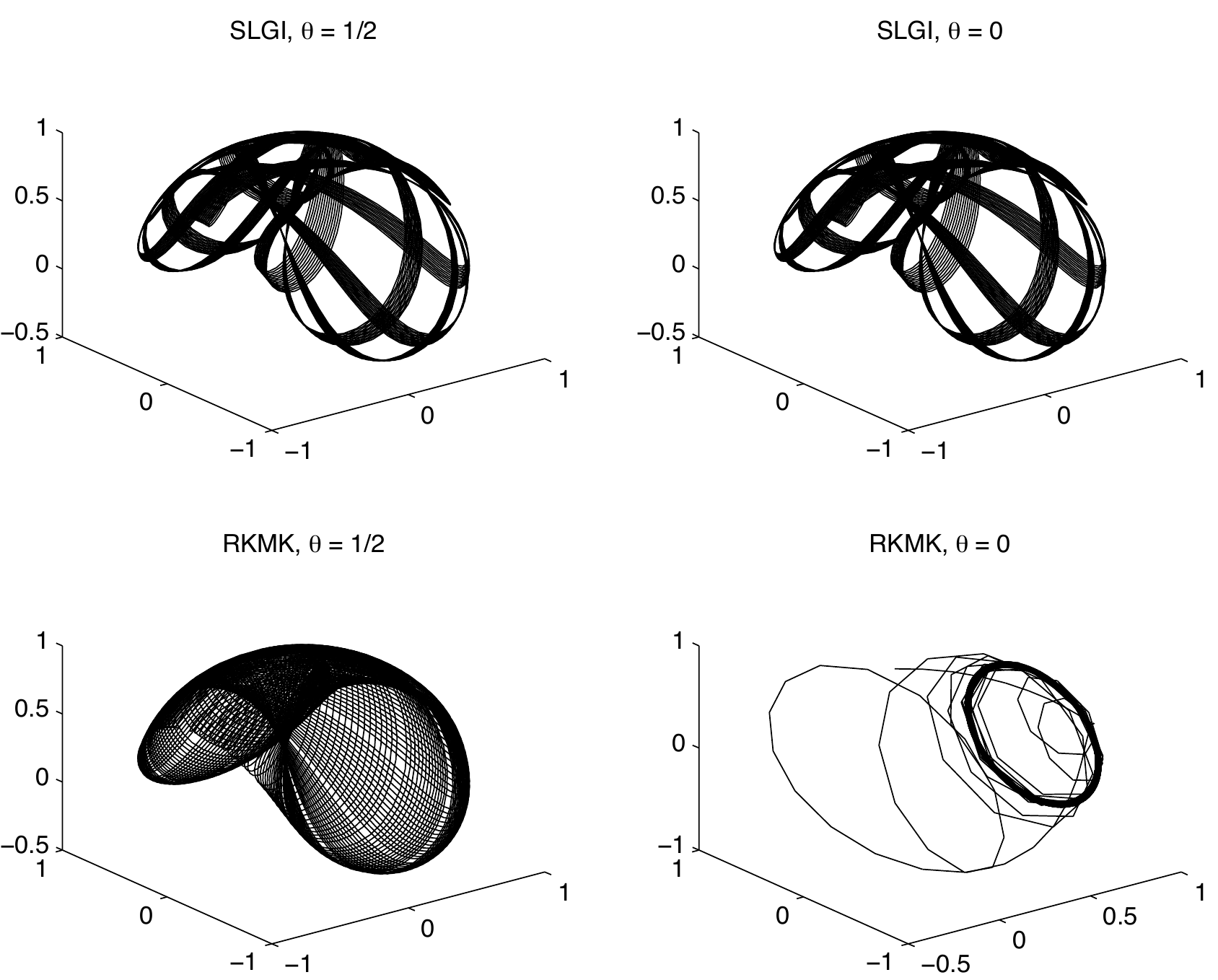}
\caption{Heavy top simulations with the symplectic (SLGI) $\theta$-methods and RKMK $\theta$-methods with $\theta=0,\tfrac{1}{2}$.
The curves show the time evolution of the centre of mass of the body.
The simulations were run over $10^5$ steps with step size $h=0.05$. See the text, page \pageref{vale:pgref}, for all other parameter values.
\label{fig:heavytop}}
\end{figure}
In Figure~\ref{fig:heavytop} we show numerical experiments for the heavy top, where the Hamiltonian is given as
\[
    H \colon \mathbf{G} \to \R, \qquad H(g,\mu) =\frac{1}{2} \langle \mu, \mathbb{I}^{-1}\mu\rangle + \mathbf{e}_3^Tgu_0,
\]
where $G = \SO(3)$, $\mathbb{I}\colon \g\rightarrow \g^*$ is the inertia tensor, here represented as a diagonal \label{vale:pgref}
$3\times 3$ matrix, $u_0$ is the initial position of the top's centre of mass, and
$\mathbf{e}_3$ is the canonical unit vector in the vertical direction. We have chosen
$\mathbb{I}=10^3\,\operatorname{diag}(1,5,6)$ and $u_0=\mathbf{e}_3$. The initial values used
were $g_0=I$ (the identity matrix), and $\mu_0=10\,\mathbb{I}\,(1,1,1)^T$.
We compare the behaviour of the symplectic schemes presented here to the Runge--Kutta--Munthe-Kaas (RKMK) method with the same coefficients.
In Figure~\ref{fig:heavytop} we have drawn the time evolution of the centre of mass,
$u_n=g_n u_0$. The characteristic band structure observed for the symplectic methods was reported in \cite{celledoni06ets}. The RKMK method with $\theta=\frac{1}{2}$ exhibits a similar behaviour, but the bands are expanding faster than for the symplectic ones.
We have also found in these experiments that none of the symplectic schemes, $\theta=0$ and $\theta=\frac{1}{2}$ have energy drift, but this is also the case for the RKMK method with $\theta=\frac{1}{2}$. This may be related to the fact that both methods are symmetric for $\theta=\frac{1}{2}$. For $\theta=0$, however, the RKMK method shows energy drift as expected. These tests were done with step size $h=0.05$ over $10^5$ steps.
See Table~\ref{tab:experiment} for a summary of the properties of the tested methods.
\begin{table}
\centering
\morespacearray
\begin{tabular}{r|c|c|c|c|}
	\cline{2-5}
	                                      & \multicolumn{2}{c|}{RKMK}     & \multicolumn{2}{c|}{SLGI}     \\
	\cline{2-5}
	                                      & $\theta = 0$ & $\theta = 1/2$ & $\theta = 0$ & $\theta = 1/2$ \\
	\hline
	\multicolumn{1}{|r|}{Symplectic}      & no           & no             & yes          & yes            \\
	\multicolumn{1}{|r|}{Symmetric}       & no           & yes            & no           & yes            \\
	\multicolumn{1}{|r|}{No energy drift} & no           & yes            & yes          & yes            \\
	\hline
\end{tabular}
\caption{Properties of the tested methods. The energy drift was observed numerically.} \label{tab:experiment}
\end{table}

\section{Discrete gradients and integral preserving methods on Lie groups}
\label{sec:discdiff}
The discrete gradient method for preserving first integrals
has to a large extent been made popular through the works of Gonzalez \cite{gonzalez96tia} and McLachlan et al.\ \cite{mclachlan99giu}. The latter proved the result that under relatively general circumstances, a differential equation which has a first integral $I(x)$ can be written in the form
\[
       \dot{x} = S(x) \nabla I(x),
\]
for some non-unique solution-dependent skew-symmetric matrix $S(x)$. The idea is to introduce a mapping which resembles the true gradient; a \emph{discrete gradient}
$\overline{\nabla}I \colon \R^d\times\R^d\rightarrow \R^d$ which is a continuous map that satisfies the following two conditions:
\begin{align*}
\overline{\nabla}I(x,x) &= \nabla I(x),\quad \forall x, \\
I(y)-I(x) &= \overline{\nabla}I(x,y)^T(y-x),\quad \forall x\neq y.
\end{align*}
An integrator which preserves $I$, that is, $I(x_n)=I(x_0)$ for all $n$ is now easily devised
as
\[
      \frac{x_{n+1}-x_n}{h} = \tilde{S}(x_n,x_{n+1})\overline{\nabla}I(x_n,x_{n+1}),
\]
where $\tilde{S}(x,y)$ is some consistent approximation to $S(x)$, i.e.\ $\tilde{S}(x,x)=S(x)$.
There exist several discrete gradients, two of the most popular are
\begin{equation} \label{eq:cavf}
     \overline{\nabla}I(x,y) = \int_0^1 \nabla I(\zeta y + (1-\zeta) x) \, \d \zeta,
\end{equation}
and
\begin{equation} \label{eq:gonz}
\overline{\nabla} I(x,y) = \nabla I\left(\frac{x+y}{2}\right) + \frac{I(y)-I(x)-\nabla I\left(\frac{x+y}{2}\right)^T(y-x)}{\lVert y-x\rVert^2} (y-x).
\end{equation}
The matrix $\tilde{S}(x,y)$ can be constructed with the purpose of increasing the accuracy of the resulting approximation, see e.g.\ \cite{quispel08anc}. 

We now generalise the concept of the discrete gradient to a Lie group $G$. We consider differential equations which can, for a given dual two-form\footnote{By dual two-form, we here mean a differential two-form on $G$ such that on each fibre of the cotangent bundle
we have $\omega_x \colon T_x^*G\times T_x^*G\rightarrow\R$, a bilinear, skew-symmetric form. Such forms are sometimes called bivectors or 2-vectors.}
 $\omega\in\Omega_2(G)$ and a function $H \colon G\rightarrow\R$ be written in the form
\begin{equation*} \label{eq:intprod}
      \dot{x} = \mathrm{i}_{\d H}\omega,
\end{equation*}
where $\mathrm{i}_{\alpha}$ is the interior product 
$\mathrm{i}_\alpha\omega(\beta)=\omega(\alpha,\beta)$ for any two one-forms
$\alpha, \beta\in\Omega^1(G)$.
The function $H$ is a first integral since
\[
    \frac{\d}{\d t} H(x(t)) = \d H_{x(t)}(\dot{x}(t))=\omega(\d H, \d H)=0.
\]
We define the \emph{trivialised discrete differential} (TDD)  of the function $H$ to be a continuous map
$\dd H \colon G\times G \rightarrow \g^*$ such that
\begin{align*}
 H(x')-H(x) &=     \langle \dd H(x,x'), \log(x' \cdot x^{-1}) \rangle,   \\
 \dd H(x,x) &= R_x^*\,\d H_x.
\end{align*}
A numerical method can now be defined in terms of the discrete differential as
\[
     x' = \exp(h \, \mathrm{i}_{\dd H(x,x')}\bar{\omega}(x,x')) \cdot x.
\]
where $\bar{\omega}$ is a continuous map from $G\times G$ into the space of exterior two-forms on $\g^*$,
$\Lambda^2(\g^*)$.
 This exterior form
is some local trivialised approximation to $\omega$, meaning that we 
impose the following  consistency condition
\begin{equation*} \label{eq:consomb}
     \bar{\omega}(x,x)(R_x^*\alpha, R_x^*\beta) = \omega_x(\alpha,\beta), \qquad \text{for all }
     \alpha, \beta\in T_x^*G.
\end{equation*}
We easily see that this method preserves $H$ exactly, since
\begin{align*}
H(x')-H(x) &= \langle\dd H(x,x'), \log(x' \cdot x^{-1})\rangle \\
                &= \langle\dd H(x,x'), h\,\mathrm{i}_{\dd H(x,x')}\bar\omega(x, x')\rangle \\
                &= h\,\bar{\omega}(x,x') (\dd H(x,x'), \dd H(x,x')) = 0.
\end{align*}
Extending \eqref{eq:cavf} to the Lie group setting, we define the following TDD:
\begin{equation*} \label{eq:tdgavf}
   \dd H(x,x') = \int_0^1 R_{\ell(\xi)}^* \, \d H_{\ell(\xi)} \, \d\xi,\qquad
   \ell(\xi) = \exp(\xi\log(x'\cdot x^{-1})) \cdot x.
\end{equation*}
Similarly, for any given inner product on $\g$, we may extend the discrete gradient \eqref{eq:gonz}
to
\begin{equation*}
\label{eq:gonz1}
    \dd H(x,x') = R_{\bar{x}}^* \, \d H_{\bar{x}} +
    \frac{H(x')-H(x)-\langle R_{\bar{x}}^* \, \d H_{\bar{x}},\eta\rangle
    }{\lVert\eta\rVert^2}\eta^\flat,\quad\eta=\log(x' \cdot x^{-1}),
    \end{equation*}
where $\bar{x}\in G$ for instance could be $\bar{x}=\exp(\eta/2) \cdot x$, a choice which would
cause $\dd H(x,x')=\dd H(x',x)$. 
The standard notation $\eta^\flat$ is used for index-lowering, the inner product $({\cdot},{\cdot})$ associates to any
element $\eta\in\g$ the dual element $\eta^\flat\in\g^*$ through $\langle\eta^\flat, \zeta\rangle = (\eta,\zeta)$, $\forall\zeta\in\g$.

Suppose that the ODE vector field $F$ is known as well as the invariant $H$.
A dual two-form $\omega$ can now be defined in terms of a Riemannian metric on 
$G$.  By index raising applied to $\d H$, we obtain the Riemannian gradient vector field $\grad H$, and we define
\[
      \omega = \frac{\grad H\wedge F}{\lVert\grad H\rVert^2}\qquad
      \Rightarrow\qquad \mathrm{i}_{\d H}\omega = F.
\]
\paragraph{Example} We consider the equations for the attitude rotation of a free rigid body expressed using Euler parameters. The set  $S^3=\{  \q\in \R^4 \mid  \lVert\q\rVert_2=1\}$ with $\q=(q_0,\mathbf{q})$, ($q_0\in\R$ and $\mathbf{q}\in\R^3$), is a Lie group  with the quaternion product 
\[
	\mathbbm{p}\cdot \q=(p_0q_0-\mathbf{p}^T\mathbf{q}, p_0\mathbf{q}+q_0\mathbf{p}+\mathbf{p}\times \mathbf{q}),
\]
with unit $\mathbbm{e}=(1, 0, 0, 0)$ and inverse
$\q_c=(q_0,\,-\mathbf{q})$. We denote by ``$\hat{\quad}$" the hat-map defined in \eqref{eq:hatmap}. 
The Lie group $S^3$ can be mapped into $\SO(3)$ by the Euler--Rodrigues map:
\[
	\mathcal{E}(\q)=I_3+2q_0\hat{\mathbf{q}}+2\hat{\mathbf{q}}^2,
\]
where $I_3$ denotes the $3\times 3$ identity matrix.
The Lie algebra $\mathfrak{s}^3$ of $S^3$ is the set of so called pure quaternions,
the elements of $\R^4$ with first component equal to zero, 
identifiable with $\R^3$ and  with $\mathfrak{so}(3)$ via the hat-map.

The equations for the attitude rotation of a free rigid body on $S^3$ read
\[
	\dot{\q}=f(\q)\cdot \q,\qquad f(\q)=\q\cdot \mathbbm{v} \cdot \q_c,
\]
and 
\[
	\mathbbm{v}=(0, \mathbf{v}),\qquad \mathbf{v} = \frac{1}{2} \mathbb{I}^{-1} \mathcal{E} (\q_c)\mathbf{m}_0,
\]
where $\mathbf{m}_0$ is the initial body angular momentum and $\mathbb{I}$ is the diagonal inertia tensor, and according to the notation previously used in this section $F(\q)=f(\q)\cdot \q$.
The energy function is
\[H(\q)=\frac{1}{2}\,\mathbf{m}_0^T\mathcal{E}(\q)\mathbb{I}^{-1}\mathcal{E}(\q_c)\mathbf{m}_0.\]

We consider the $\R^3$ Euclidean inner product  as metric in the Lie algebra $\mathfrak{s}^3$, and obtain by right translation a Riemannian metric on $S^3$.  The Riemannian gradient  of $H$ with respect to this metric is then
\[\grad H=(I_4-\q\q^T)\nabla H,\]
where $I_4$ is the identity in $\R^{4\times4}$ and $\nabla H$ is the usual gradient of $H$ as a function from $\R^4$ to $\R$. 
We identify $\mathfrak{s}^3$ with its dual, and using $\grad H$  in \eqref{eq:gonz} we obtain the (dual) discrete differential $\overline{\grad H}(\q,\q')\in \mathfrak{s}^3$.

The 
two-form
$\omega = \frac{\grad H\wedge F}{\lVert\grad H\rVert^2}$ with respect to the right trivialisation can be identified with the $4\times 4$ skew-symmetric matrix 
\[\omega_R(\q)=\frac{\xi\,\gamma^T-\gamma\,\xi^T}{\lVert\gamma\rVert^2},\qquad \xi ,\, \gamma\in\mathfrak{s}^3,\qquad \xi\cdot \q=F(\q), \qquad \gamma \cdot \q=\grad H(\q),\]
where $\omega_R(\q)$ has first row and first column equal to zero. 
We choose $\bar{\omega}$ to be
\[\bar{\omega}(\q,\q')=\omega_R(\bar{\q}),\qquad \bar{\q}=\exp(\eta/2) \cdot \q,\qquad \eta=\log (\q'\cdot \q_c),\]
i.e.\ $\omega_R$ frozen at the mid-point $\bar{\q}$.
 The energy-preserving Lie group method of second order is
 \[\q'=\exp( h\,\bar{\omega}(\q,\q') \overline{\grad H}(\q,\q'))\cdot \q,\]
 and $\exp$ is the exponential map from $\mathfrak{s}^3$ to $S^3$ with $\log \colon S^3\rightarrow \mathfrak{s}^3$ as its inverse, defined locally around the identity.
 
\begin{figure}\centering
\includegraphics[width=\textwidth]{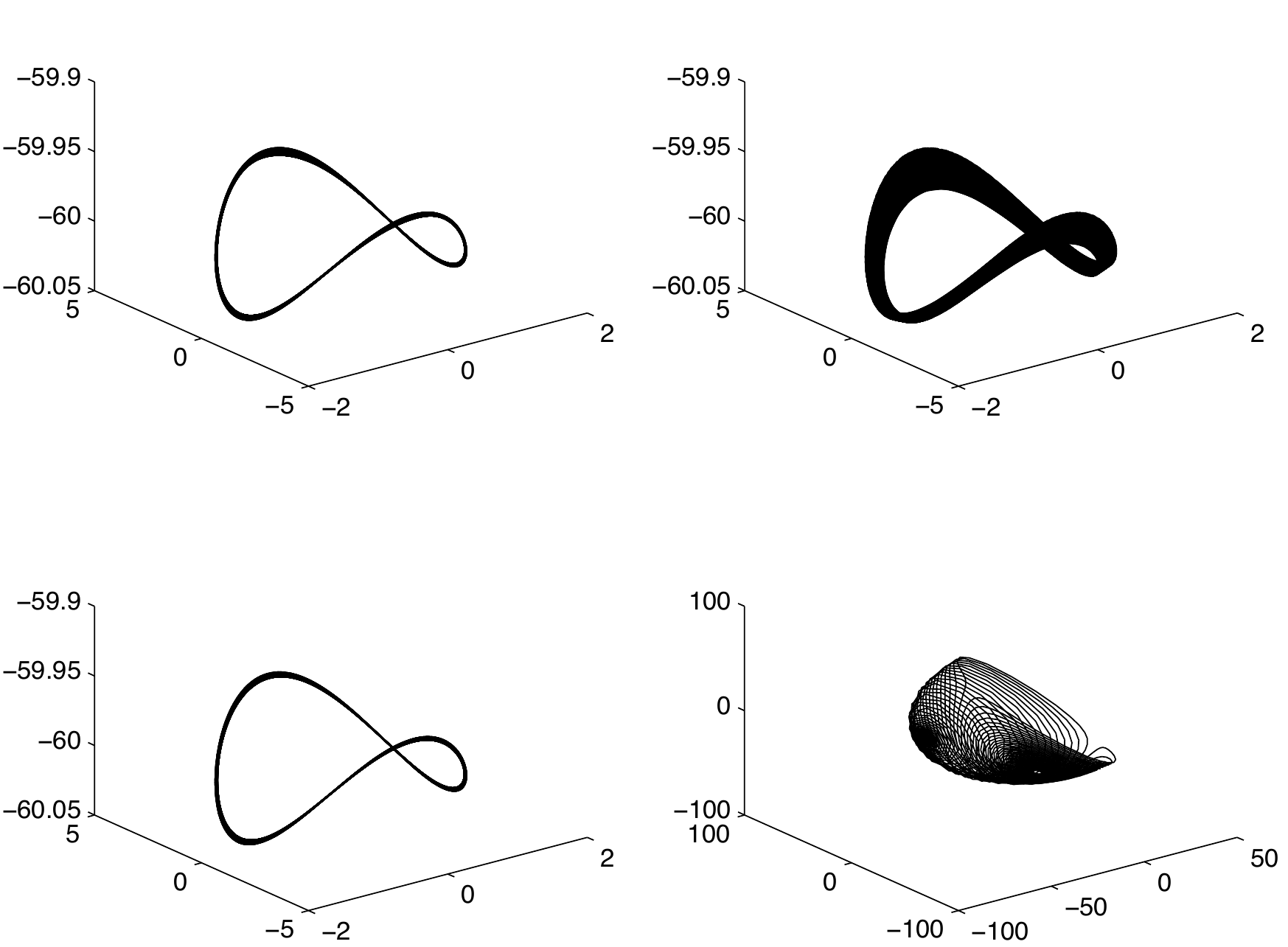}
\caption{Free rigid body angular momentum, time interval $[0,1000]$,  moments of inertia $I_1=1$, $I_2=5$, $I_3=60$, initial angular velocity $\mathbb{I}\,\mathbf{m}_0=(1, 1/2, -1)^T$. (Top left) energy-preserving Lie group method, $h=1/64$; (top right) \texttt{ode45} with tolerances $10^{-6}$;  (bottom left) \texttt{ode45} with tolerances $10^{-14}$; (bottom right) Heun RKMK, $h=1/64$.\label{fig:EP}}
\end{figure}

 In Figure~\ref{fig:EP} we plot the body angular momentum vector $\mathbf{m}=\mathcal{E}(\q_c)\mathbf{m}_0$ on a time interval $[0,T]$, $T=1000$, for four different methods: the Lie group energy-preserving integrator just described (top left),  the built-in \textsc{Matlab} routine \texttt{ode45} with absolute and relative tolerance $10^{-6}$ (top right); the \texttt{ode45} routine with tolerances $10^{-14}$ (bottom left); and the explicit Heun RKMK Lie group method (bottom right).
The two Lie group methods both have order $2$. The energy preserving method is both symmetric, energy preserving and it preserves the constraint $\lVert\q\rVert_2=1$. The Lie group integrators use a step-size $h=1/64.$ The solution of the built-in \textsc{Matlab} routine at high precision is qualitatively similar to the highly accurate solution  produced by \textsc{Matlab} with tolerances $10^{-14}$. The energy error is also comparable for these two experiments. 
The performance of other \textsc{Matlab} built-in routines we tried was worse than for \texttt{ode45}. We remark that the equations are formulated as differential equations on $S^3$, a formulation of the problem in form of a differential algebraic equation would possibly have improved the performance of the \textsc{Matlab} built-in routines. However it seems that the preservation of the constraint alone can not guarantee the good performance of the method. In fact  the explicit (non-symmetric) Lie group integrator preserves the constraint $\lVert\q\rVert_2=1$, but performs poorly on this problem (see Figure~\ref{fig:EP} bottom right). The cost per step of the explicit Lie group integrator is much lower than for the energy-preserving symmetric Lie group integrator.

\smallskip

We have given an introduction to Lie group integrators for differential equations on manifolds using the notions of frames and Lie group actions. A few application areas have been discussed. An interesting observation is that when the Lie group featuring in the method can be chosen to be the Euclidean group, the resulting integrator always reduce to some well-known numerical scheme like for instance a Runge-Kutta method. In this way, one may think of the Lie group integrators as a superset of the traditional integrators and a natural question to ask is whether the Euclidean choice will always be superior to any other Lie group and group action.

Lie group integrators that are symplectic for Hamiltonian problems in the general setting presented here are, as far as we know, not known. However, we have shown that such methods exist in the important case of Hamiltonian problems on the cotangent bundle of a Lie group. There are however still many open questions regarding this type of schemes, like for instance how to obtain high order methods.

The preservation of first integrals in Lie group integrators has been achieved in the literature by imposing a group action in which the orbit, i.e. the reachable set of points, is contained in the level set of the invariant. But it is not always convenient to impose such group actions, and we have here suggested a type of Lie group integrator which can preserve any prescribed invariant for the case where the manifold is a Lie group acting on itself by left or right multiplication. An interesting idea to pursue is the generalisation of this approach to arbitrary smooth manifolds with a group action.

\section*{Acknowledgements}
This research was supported by a Marie Curie International Research Staff Exchange Scheme Fellowship within the 7th European Community Framework
Programme.
The authors would like to acknowledge the support from the GeNuIn Applications and SpadeAce projects funded by the Research Council of Norway,  and most of the work was carried out while the authors were visiting  Massey University, Palmerston North, New Zealand and La Trobe University, Melbourne, Australia.

\bibliographystyle{elsarticle-harv}
\bibliography{geom_int,brynbib,elena}

\end{document}